\documentclass[11pt]{article}
\usepackage{amsmath,amsthm,amsfonts,amssymb,amscd, amsxtra, mathrsfs,textcomp,verbatim,inputenc}
\usepackage{url}
\usepackage{xcolor}
\usepackage{mathtools}
\usepackage{todonotes} 
\usepackage{cancel}
\usepackage{booktabs}
\usepackage[shortlabels]{enumitem}
\usepackage{soul} 
\usepackage[margin=2.5 cm,nohead]{geometry}

\usepackage{algorithm}
\theoremstyle{plain}
\theoremstyle{definition}
\newtheorem{theorem}{Theorem}
\newtheorem{lemma}[theorem]{Lemma}
\newtheorem{definition}[theorem]{Definition}
\newtheorem{corollary}[theorem]{Corollary}
\newtheorem{proposition}[theorem]{Proposition}

\newtheorem{remark}{Remark}

\newtheorem{assumption}[theorem]{Assumption}

\newtheorem*{properties*}{Properties and Examples of B-subdifferential}
\begin{document}
\title{Inexact Proximal Method on Hadamard Manifolds}
\author{
M. S. Louzeiro\thanks{
Institute of Mathematics and Statistics, Federal University of Goias,
Avenida Esperan\c{c}a, s/n, Campus II, Goi\^ania, GO - 74690-900, Brazil
(E-mail: {\tt mauriciolouzeiro@ufg.br})
}
\and
R. Serra\thanks{
Department of Mathematics, Federal University of Piau\'i,
Teresina, PI 64049-550, Brazil
(E-mail: {\tt rayserra@ufpi.edu.br}, {\tt sandoel.vieira@ufpi.edu.br})
}
\and
S. Vieira\footnotemark[2]
}
\maketitle

\begin{abstract}
We propose and analyze two inexact proximal point methods for unconstrained
geodesically convex optimization on Hadamard manifolds.
The first, a standard inexact proximal point method (SIPPM), accepts
approximate subproblem solutions directly as new iterates and applies to
arbitrary Hadamard manifolds under a summability condition on the relative
error tolerances. The second, a hybrid inexact proximal point method
(IHPPM), augments each inexact proximal step with a metric projection onto
a half-space-like separating set, allowing the summability requirement to
be relaxed to a fixed relative error tolerance at the cost of assuming
geodesic convexity of the separating sets. For both methods we establish
full asymptotic convergence, global $\mathcal{O}(\varepsilon^{-2})$
iteration-complexity bounds for approximate stationarity, and, under a
local growth condition, local $Q$-linear convergence rates---with
$Q$-superlinear acceleration for SIPPM when the proximal parameters tend
to zero and a refined contraction factor for IHPPM. We further prove that
the geodesic convexity of the separating sets holds automatically on
Hadamard surfaces and on Hadamard manifolds of constant sectional
curvature.
\end{abstract}




\section{Introduction}

Let \(\mathbb{M}\) be a finite-dimensional Hadamard manifold equipped with a Riemannian metric \(\langle \cdot, \cdot \rangle\), and consider the unconstrained convex optimization problem
\begin{equation}
    \min_{p\in\mathbb{M}} f(p),
    \label{eq:problem}
\end{equation}
where \(f: \mathbb{M} \to \mathbb{R}\) is a geodesically convex function. Throughout this paper, we assume that the solution set of
problem~\eqref{eq:problem} is nonempty. Optimization on Riemannian manifolds has attracted considerable attention over the past decades, evolving into a mature framework for handling problems with intrinsic geometric constraints. Rather than enforcing ambient non-convex constraints, manifold optimization exploits the underlying geometry to convert constrained Euclidean problems into unconstrained Riemannian ones \cite{absil2008optimization, boumal2020introduction,Udriste1994}.

In particular, optimization on non-positively curved spaces—such as Hadamard manifolds and \(\mathrm{CAT}(0)\) metric spaces—has found crucial applications across diverse modern domains. Notable examples include metric learning and low-rank matrix completion \cite{vandereycken2013low}, computer vision and medical image analysis via symmetric positive definite (SPD) matrices \cite{pennec2006riemannian,cherian2016riemannian}, principal geodesic analysis and geometric statistics \cite{fletcher2004principal}, and optimal transport or dictionary learning on hyperbolic spaces \cite{nickel2017poincare, sala2018representation}. In these contexts, Hadamard manifolds offer a remarkably rich framework: complete geodesic connectivity, the global diffeomorphic nature of the exponential mapping, and the non-expansiveness of geodesics allow fundamental Euclidean tools to be extended with robust convergence properties.

Among the foundational algorithms for solving \eqref{eq:problem}, the Proximal Point Method (PPM) plays a central role due to its exceptional stability and regularizing properties. Introduced in the Riemannian setting by Ferreira and Oliveira \cite{Ferreira2002} (see also \cite{bacak2014convex}), the exact proximal scheme updates the current iterate \(p_k\) by solving the subproblem
\begin{equation}\label{eq:intro-exact-prox}
    p_{k+1} = \arg\min_{p\in\mathbb{M}} \left\{ f(p) + \frac{\mu_k}{2} d^2(p, p_k) \right\},
\end{equation}
where \(\mu_k > 0\) is a regularization parameter and \(d(\cdot,\cdot)\) denotes the Riemannian distance. The optimality condition for \eqref{eq:intro-exact-prox} corresponds to the proximal inclusion \(\mu_k \log_{p_{k+1}} p_k \in \partial f(p_{k+1})\). Although exact PPM generates sequence iterates with strong Fej\'er monotonicity properties toward the solution set \(S^* := \arg\min_{\mathbb{M}} f\), calculating \(p_{k+1}\) exactly at each iteration is often computationally prohibitive when no closed-form expression for the proximal operator is available. This computational bottleneck is particularly pronounced in high-dimensional signal processing, deep hyperbolic representation learning, or large-scale SPD covariance estimation. This motivates the design of \emph{inexact} proximal point strategies, where the subproblems are solved only up to a controlled error tolerance using inner optimization routines.

In this paper, we propose and analyze two complementary inexact proximal point algorithms for solving \eqref{eq:problem} on Hadamard manifolds, providing a unified perspective on the trade-offs between computational accuracy, geometric restrictions, and convergence guarantees:

\begin{enumerate}[label=(\roman*)]
    \item \textbf{Standard Inexact Proximal Point Method (SIPPM):} To avoid geometric restrictions on the separating sets altogether, we first investigate a direct inexact scheme where the approximate subproblem solution itself becomes the next iterate, i.e., \(p_{k+1} = q_k\). At iteration \(k\), an approximate proximal triple \((q_k, v_k, \varepsilon_k)\) is computed satisfying
    \[
        0 = v_k - \mu_k \log_{q_k} p_k + \varepsilon_k, \;\; v_k \in \partial f(q_k),
    \]
    subject to an iteration-dependent relative error bound \(\|\varepsilon_k\| \le \sigma_k \mu_k d(q_k, p_k)\). To prevent error accumulation and ensure that sequence iterates maintain Fej\'er monotonicity toward \(S^*\), the sequence of relative error tolerances \(\{\sigma_k\}_{k\in\mathbb{N}} \subset [0, 1)\) cannot be chosen arbitrarily. It must be dynamic and satisfy the summability condition \(\sum_{k=0}^{\infty} \sigma_k < +\infty\) \cite{Rockafellar76}. The primary advantage of SIPPM is its universal applicability: it operates on \emph{arbitrary} finite-dimensional Hadamard manifolds without any curvature restrictions beyond nonpositivity.

    \item \textbf{Inexact Hybrid Proximal Point Method (IHPPM):} Inspired by the seminal work of Solodov and Svaiter \cite{SolodovBenar1999} in Hilbert spaces, this algorithm combines an inexact proximal step with a metric projection onto a separating set containing \(S^*\). A fundamental advantage of IHPPM is that it relaxes the stringent requirement on the error parameters: instead of a vanishing summable sequence, one can use a \emph{fixed} tolerance parameter \(\sigma_k \equiv \sigma \in [0, 1)\) for all \(k \in \mathbb{N}\). The current iterate \(p_k\) is then projected onto the half-space-like separating set
    \[
        \mathcal{S}_{q_k, v_k} := \{ p \in \mathbb{M} : \langle v_k, \log_{q_k} p \rangle \le 0 \}.
    \]
 The metric projection step acts as a geometric safeguard, absorbing subproblem inexactness and restoring full Fej\'er monotonicity with respect to \(S^*\). As observed in \cite{GarcigaIusemSvaiter2001}, without such a projection step, inexact proximal iterations using a fixed relative error $\sigma > 0$ can generate unbounded sequences even in Euclidean spaces. Furthermore, when \(\mathbb{M}\) is the real hyperbolic space \(\mathbb{H}^n\), we derive an explicit closed-form formula for this metric projection, greatly facilitating its practical implementation.
 
 While \(\mathcal{S}_{q,v}\) is not automatically geodesically convex on general Hadamard manifolds \cite{wang2016some}, our framework and the hybrid method proposed in \cite{DFLJ2026} offer complementary perspectives on navigating this geometric property. While the approach in \cite{DFLJ2026} operates on arbitrary Hadamard manifolds by leveraging a specific structural condition on \(f\), our hybrid scheme guarantees applicability to the full class of  geodesically convex functions \(f\). To accommodate arbitrary convex objectives without structural restrictions, we address the underlying geometry directly by identifying conditions under which \(\mathcal{S}_{q_k,v_k}\) is geodesically convex. In particular, we establish global convergence when \(\mathbb{M}\) features constant sectional curvature (such as hyperbolic spaces) or, more generally, when the separating boundary \(\mathcal H_{q_k,v_k} := \exp_{q_k}(v_k^\perp)\) forms a totally geodesic hypersurface. Notably, this geometric setting is naturally realized on any two-dimensional Hadamard manifold (Hadamard surfaces), enabling our hybrid method to seamlessly accommodate broad classes of spaces with variable nonpositive curvature for arbitrary convex objective functions.
\end{enumerate}

For both methods, we prove that under the standard regularization condition \(\sum_{k=0}^{\infty} \mu_k^{-2} = +\infty\), the generated sequences converge globally to a point in \(S^*\). Furthermore, we derive iteration-complexity bounds for finding \(\epsilon\)-stationary points as well as functional value residual bounds. Under a local growth condition on \(f\), we establish Q-linear convergence rates for both methods and Q-superlinear convergence for the direct method when \(\mu_k \to 0\). For the hybrid method, we exploit a detailed case analysis of the relative error criterion to establish a refined local linear convergence rate estimate.

The remainder of this paper is organized as follows. Section~\ref{sec:int.r} reviews fundamental geometric properties of Hadamard manifolds, metric projections, Riemannian subdifferentials, and Fej\'er convergence properties. Section~\ref{app:direct-inexact-proximal} is devoted to the SIPPM, presenting its global convergence analysis, iteration complexity estimates, and local rate of convergence under the error bound condition. Section~\ref{sec:proxLine} introduces the IHPPM, detailing its global convergence theory, complexity bounds, local rate analysis, and a comprehensive discussion on the geodesic convexity of the separating sets \(\mathcal{S}_{q,v}\). Section~\ref{sec:conclusions} provides concluding remarks and discusses future directions. Finally, two self-contained appendixes complement our theoretical framework: Appendix~\ref{sec:hyperbolic-projection} derives an explicit closed-form formula for the metric projection onto a hyperbolic hyperplane in the hyperboloid model of \(\mathbb{H}^n\), while Appendix~\ref{App:ass:conv.su.} establishes precise geometric characterizations and sufficient conditions for the geodesic convexity of the separating sets \(\mathcal{S}_{q,v}\) on Hadamard surfaces and manifolds of constant sectional curvature.

\section{Preliminaries}\label{sec:int.r}

In this section, we recall the notation and basic facts on Hadamard manifolds that will be used throughout the paper; for further details, see, e.g., \cite{bacak2014convex,Manfredo1992,Sakai1996,Udriste1994}.

\subsection{Hadamard manifolds}

Throughout the paper, $\mathbb{M}$ denotes a finite-dimensional \emph{Hadamard manifold}, that is, a complete, simply connected Riemannian manifold with nonpositive sectional curvature. We denote by $T_p\mathbb{M}$ the tangent space to $\mathbb{M}$ at $p$, by $T\mathbb{M}:=\bigcup_{p\in\mathbb{M}}\{p\}\times T_p\mathbb{M}$ the tangent bundle, and by $\langle\cdot,\cdot\rangle$ and $\|\cdot\|$ the Riemannian metric and its associated norm, respectively. The length of a piecewise smooth curve $\gamma$ is denoted by $\ell(\gamma)$, and $d$ denotes the Riemannian distance. The metric space $(\mathbb{M},d)$ is complete and induces the manifold topology.
By the Cartan--Hadamard theorem, for every $p\in\mathbb{M}$ the exponential map $\exp_p\colon T_p\mathbb{M}\to\mathbb{M}$ is a diffeomorphism. Its inverse is denoted by $\log_p\colon\mathbb{M}\to T_p\mathbb{M}$. For any $p,q\in\mathbb{M}$, there exists a unique geodesic segment joining $p$ to $q$, given by
$\gamma_{pq}(t):=\exp_p\bigl(t\log_p q\bigr)$, $t\in[0,1]$,
and it satisfies $\ell(\gamma_{pq})=d(p,q)=\|\log_p q\|$. When no confusion can arise, we write $\log_p q$ instead of $\log_p(q)$.

The next result is the well-known comparison theorem for geodesic triangles in Hadamard manifolds, see \cite[Chapter~4, Proposition~4.5]{Sakai1996}.

\begin{lemma}\label{le:CosLawF}
For all $x,y,z\in\mathbb{M}$,
\[
\|\log_x y-\log_x z\|^2
= d^2(x,y)+d^2(x,z)-2\langle\log_x y,\log_x z\rangle
\le d^2(y,z).
\]
\end{lemma}

\subsection{Convex sets and metric projections}

A subset $C\subset\mathbb{M}$ is \emph{convex} if, for all $p,q\in C$, the geodesic segment $\gamma_{pq}$ is contained in $C$. If $C$ is closed and convex, then every $q\in\mathbb{M}$ admits a unique nearest point in $C$, called the \emph{metric projection} of $q$ onto $C$ and denoted by
$
\mathcal{P}_C(q):=\operatorname*{arg\,min}_{x\in C}\, d(q,x).
$
The following proposition characterizes the projection, its proof can be found in \cite{Ferreira2002}.

\begin{proposition}\label{prop:proj-ineq}
Let $C\subset\mathbb{M}$ be a closed convex set and $q\in\mathbb{M}$. Then, for every $p\in C$,
\[
\langle\log_{\mathcal{P}_C(q)}q,\ \log_{\mathcal{P}_C(q)}p\rangle\le 0.
\]
\end{proposition}

\subsection{Convex functions and subgradients}

A function $f\colon\mathbb{M}\to\mathbb{R}$ is \emph{convex} if
\[
f(\gamma_{pq}(t))\le(1-t)f(p)+t f(q),\;\;\forall\,p,q\in\mathbb{M},\ \forall\,t\in[0,1].
\]

\begin{definition}\label{def:subgradient}
A vector $s\in T_q\mathbb{M}$ is a \emph{subgradient} of $f\colon\mathbb{M}\to\mathbb{R}$ at $q$ if
$f(p)\ge f(q)+\langle s,\log_q p\rangle$ for all $p\in\mathbb{M}$.
The set of all subgradients of $f$ at $q$, denoted by $\partial f(q)$, is called the \emph{subdifferential} of $f$ at $q$.
\end{definition}

It is well known that a function $f\colon\mathbb{M}\to\mathbb{R}$ is convex if and only if $\partial f(q)\neq\emptyset$ for every $q\in\mathbb{M}$ and the subgradient inequality in Definition~\ref{def:subgradient} holds for all $p,q\in\mathbb{M}$ and all $s\in\partial f(q)$.

\begin{definition}\label{def:epsilon-critical-point}
Let $\epsilon\ge 0$. A point $p\in\mathbb{M}$ is an \emph{$\epsilon$-critical point} of a convex function $f$ if there exists $v\in\partial f(p)$ such that $\|v\|\le\epsilon$.
\end{definition}

The following result is an immediate consequence of \cite[Proposition~2.5]{WangLiWangYao2015}.

\begin{proposition}\label{prop:subdifferential-continuity}
Let $f\colon\mathbb{M}\to\mathbb{R}$ be convex and let $\{p_k\}_{k\in\mathbb{N}}\subset\mathbb{M}$ be a sequence converging to $ p\in\mathbb{M}$. If $s_k\in\partial f(p_k)$ for every $k\in\mathbb{N}$, then $\{s_k\}_{k\in\mathbb{N}}$ is bounded and each of its cluster points belongs to $\partial f( p)$.
\end{proposition}

\subsection{Strong convexity}

A function $f\colon\mathbb{M}\to\mathbb{R}$ is \emph{$m$-strongly convex}, with $m>0$, if
\[
f(\gamma_{pq}(t))\le(1-t)f(p)+t f(q)-\frac{m}{2}\,t(1-t)\,d^2(p,q),\;\;\forall\,p,q\in\mathbb{M},\ \forall\,t\in[0,1].
\]
Equivalently, $f$ is $m$-strongly convex if and only if
\begin{equation}\label{eq:strong-subgradient-inequality}
f(p)\ge f(q)+\langle v,\log_q p\rangle+\frac{m}{2}\,d^2(p,q),\;\;\forall\,p,q\in\mathbb{M},\ \forall\,v\in\partial f(q).
\end{equation}
The proof of this equivalence follows by standard arguments, analogous to those in \cite[Theorem~4.10, p.~76]{Udriste1994}, and is therefore omitted.

\begin{lemma}\label{lem:cat0}
Let $f\colon\mathbb{M}\to\mathbb{R}$ be convex, $p\in\mathbb{M}$, and $m>0$. Then the function
\[
\phi(q):=f(q)+\frac{m}{2}\,d^2(q,p)
\]
is $m$-strongly convex. Consequently, $\phi$ has a unique minimizer $x\in\mathbb{M}$, and
\begin{equation}\label{eq:strong-min}
\phi(q)\ge\phi(x)+\frac{m}{2}\,d^2(q,x),\;\;\forall\,q\in\mathbb{M}.
\end{equation}
\end{lemma}

\begin{proof}
Since $\mathbb{M}$ has nonpositive sectional curvature, the function $q\mapsto\frac12 d^2(q,p)$ is $1$-strongly convex. As $f$ is convex, $\phi$ is $m$-strongly convex. Being convex and finite-valued, $\phi$ is continuous. Moreover, it is coercive. Hence $\phi$ attains its minimum, and the minimizer $x$ is unique by strict convexity. Since $\phi(q)\ge\phi(x)=\phi(x)+\langle 0,\log_x q\rangle$ for all $q\in\mathbb{M}$, we have $0\in\partial\phi(x)$. Applying \eqref{eq:strong-subgradient-inequality} to $\phi$ at the point $x$ with $v=0$ yields \eqref{eq:strong-min}.
\end{proof}

\subsection{Fej\'er convergence}

\begin{definition}\label{def:fejer}
Let $(X,d)$ be a complete metric space. A sequence $\{y_k\}_{k\ge 0}\subset X$ is \emph{Fej\'er convergent} to a set $W\subset X$ if
\[
d^2(y_{k+1},w)\le d^2(y_k,w),\;\;\forall\,w\in W,\ \forall\,k\ge 0.
\]
\end{definition}

\begin{proposition}\label{prop:fejer}
Let $(X,d)$ be a complete metric space and let $\{y_k\}\subset X$ be Fej\'er convergent to a nonempty set $W\subset X$. Then $\{y_k\}$ is bounded. Moreover, if a cluster point $y$ of $\{y_k\}$ belongs to $W$, then $\lim_{k\to\infty}y_k=y$.
\end{proposition}
\section{Standard Inexact Proximal Point Method}
\label{app:direct-inexact-proximal}

Before presenting our approach, we recall the classical exact proximal point
method on Hadamard manifolds, introduced in~\cite{Ferreira2002}.

\begin{algorithm}[H]
\begin{footnotesize}
\begin{description}
\item[Step 0.]
Choose a sequence of positive parameters \(\{\mu_k\}_{k\in\mathbb N}\subset(0,+\infty)\)
and an initial point \(p_0\in\mathbb M\). Set \(k\gets0\).
\item[Step 1.]
Given \(p_k\in\mathbb M\), compute \(p_{k+1}\in\mathbb M\) satisfying
\begin{equation}
\label{eq:prox-exact}
\mu_k\log_{p_{k+1}}p_k\in\partial f(p_{k+1}).
\end{equation}
\item[Step 2.]
If \(p_{k+1}=p_k\), \textbf{stop} and return \(p_k\); otherwise, set
\(k\gets k+1\) and return to \textbf{Step~1}.
\end{description}
\caption{Exact Proximal Point Method (EPPM)}
\label{Alg:ExactPPM}
\end{footnotesize}
\end{algorithm}

Since the Riemannian gradient of \((1/2) d^2(\cdot,p_k)\) at \(x\) is \(-\log_x p_k\), inclusion~\eqref{eq:prox-exact} corresponds to the optimality condition for the proximal subproblem
\[
p_{k+1}=\operatorname*{arg\,min}_{p\in\mathbb M}
\Bigl\{f(p)+\frac{\mu_k}{2}\,d^2(p,p_k)\Bigr\}.
\]
Its objective function is strongly convex, ensuring that \(p_{k+1}\) is uniquely determined. Although this exact formulation provides a natural framework for convergence analysis, computing \(p_{k+1}\) exactly can be computationally demanding in practice: when the subproblem lacks a closed-form solution, one must resort to an iterative procedure, rendering the exact enforcement of~\eqref{eq:prox-exact} unnecessarily restrictive.
This motivates an inexact variant in which~\eqref{eq:prox-exact} is satisfied
only up to a controlled error. Instead of the exact proximal point, we seek a
triple \((q_k,v_k,\varepsilon_k)\) with \(v_k\in\partial f(q_k)\) and
\(\varepsilon_k\in T_{q_k}\mathbb M\) such that
$
0=v_k-\mu_k\log_{q_k}p_k+\varepsilon_k .
$
The vector \(\varepsilon_k\) is the residual of the proximal inclusion; if
\(\varepsilon_k=0\), the optimality condition of Algorithm~\ref{Alg:ExactPPM}
is recovered. The central question is how much inexactness can be tolerated at
each iteration without losing the convergence properties of the method.
Following the relative-error framework of Rockafellar~\cite{Rockafellar76}, we
control the residual by a relative criterion with tolerances
\(\{\sigma_k\}_{k\in\mathbb N}\) that satisfy the summability condition
formalized in Assumption~\ref{ass:sigma-summable} below. The resulting scheme,
stated in Algorithm~\ref{Alg:SIPPM}, can be implemented and analyzed on an
arbitrary Hadamard manifold.

\begin{algorithm}[H]
\begin{footnotesize}
\begin{description}
\item[Step 0.]
Choose \(p_0\in\mathbb M\), a constant \(\bar\sigma\in[0,1)\), and a
sequence \(\{\sigma_k\}_{k\in\mathbb N}\subset[0,\bar\sigma]\).
Set \(k\gets0\).
\item[Step 1.]
Choose \(\mu_k>0\) and compute a triple \((q_k,v_k,\varepsilon_k)\), with
\(v_k,\varepsilon_k\in T_{q_k}\mathbb M\), satisfying
\begin{equation}
\label{eq:direct-proximal-inclusion}
0=v_k-\mu_k\log_{q_k}p_k+\varepsilon_k,
\;\;
v_k\in\partial f(q_k),
\end{equation}
and
\begin{equation}
\label{eq:rockafellar-relative-error}
\|\varepsilon_k\|\leq\sigma_k\,\mu_k\,d(q_k,p_k).
\end{equation}
\textbf{Stopping test:} if \(v_k=0\) or \(q_k=p_k\), \textbf{stop} and return \(q_k\).
\item[Step 2.]
Set \(p_{k+1}:=q_k\) and \(k\gets k+1\), and return to \textbf{Step~1}.
\end{description}
\caption{Standard Inexact Proximal Point Method (SIPPM)}
\label{Alg:SIPPM}
\end{footnotesize}
\end{algorithm}

The stopping test is sound, since it can only be triggered at a solution of
problem~\eqref{eq:problem}. Indeed, if \(v_k=0\), then \(0\in\partial f(q_k)\),
so \(q_k\) solves~\eqref{eq:problem}. If instead \(q_k=p_k\), then
\eqref{eq:rockafellar-relative-error} gives \(\varepsilon_k=0\), and
\eqref{eq:direct-proximal-inclusion} yields \(v_k=0\). Hence \(q_k\) is again a
solution.

\begin{assumption}
\label{ass:sigma-summable}
The sequence \(\{\sigma_k\}_{k\in\mathbb N}\) satisfies
\[
\Sigma:=\sum_{k=0}^{+\infty}\sigma_k<+\infty,
\;\;
0\leq\sigma_k\leq\bar\sigma<1 .
\]
\end{assumption}

\begin{remark}
\label{rem:sigma-sq-summable}
Since \(0\leq\sigma_k\leq\bar\sigma<1\), Assumption~\ref{ass:sigma-summable}
implies \(\sigma_k^2\leq\bar\sigma\,\sigma_k\), and hence
\[
\Sigma_2:=\sum_{k=0}^{+\infty}\sigma_k^2\leq\bar\sigma\,\Sigma<+\infty .
\]
\end{remark}

\begin{remark}
\label{rem:fixed-error-failure}
Assumption~\ref{ass:sigma-summable} is in line with the classical treatment of
computational errors in inexact proximal point methods for maximal monotone
operators~\cite{Rockafellar76}. In particular, it forces \(\sigma_k\to0\), so
the admissible relative error becomes increasingly stringent along the
iterations. This requirement is not an artifact of the analysis: Solodov and
Svaiter~\cite[p.~4]{SolodovBenar1999} and Garc\'ia Otero, Iusem, and
Svaiter~\cite{GarcigaIusemSvaiter2001} exhibit counterexamples, even for
operators of the form \(\partial f\), in which the direct iteration
\(p_{k+1}=q_k\) with a \emph{fixed} relative error \(\sigma>0\) fails to
converge. In Section~\ref{sec:proxLine}, we show that an additional projection
onto a separating half-space removes this limitation and allows a constant
tolerance \(\sigma\in[0,1)\).
\end{remark}

\subsection{Asymptotic Convergence Analysis for Algorithm~\ref{Alg:SIPPM}}

We start with the basic quasi-Fej\'er estimates for the direct update
\(p_{k+1}=q_k\).

\begin{lemma}
\label{le:direct-fejer-bounds}
Let \(\{p_k\}_{k\in\mathbb N}\) be a sequence generated by
Algorithm~\ref{Alg:SIPPM}, let \(p^*\in S^*\), and suppose that
Assumption~\ref{ass:sigma-summable} holds. Then, for every \(k\in\mathbb N\),
\begin{enumerate}
\item[(i)] \(d(p_k,p^*)\leq d(p_0,p^*)\,e^{\frac{\Sigma_2}{2(1-\bar\sigma^2)}}\);
\item[(ii)] \(d^2(p_{k+1},p_k)\leq d^2(p_k,p^*)-d^2(p_{k+1},p^*)+4\sigma_k\,d^2(p_0,p^*)\,e^{\frac{\Sigma_2}{1-\bar\sigma^2}}\);
\item[(iii)] \(\|v_k\|\leq(1+\bar\sigma)\,\mu_k\,d(p_{k+1},p_k)\).
\end{enumerate}
\end{lemma}

\begin{proof}
(i) Since \(p_{k+1}=q_k\) and \(p^*\in S^*\), relation~\eqref{eq:direct-proximal-inclusion}
and Definition~\ref{def:subgradient} yield
\[
\bigl\langle\mu_k\log_{p_{k+1}}p_k-\varepsilon_k,\ \log_{p_{k+1}}p^*\bigr\rangle
=\langle v_k,\log_{p_{k+1}}p^*\rangle
\leq f(p^*)-f(p_{k+1})\leq0 .
\]
Rearranging and applying the Cauchy--Schwarz inequality together with the
relative-error criterion~\eqref{eq:rockafellar-relative-error}, we obtain
\begin{equation}
\label{eq:direct-inner-product-estimate}
\bigl\langle\log_{p_{k+1}}p_k,\ \log_{p_{k+1}}p^*\bigr\rangle
\leq\frac{1}{\mu_k}\langle\varepsilon_k,\log_{p_{k+1}}p^*\rangle
\leq\frac{1}{\mu_k}\|\varepsilon_k\|\,d(p_{k+1},p^*)
\leq\sigma_k\,d(p_{k+1},p_k)\,d(p_{k+1},p^*).
\end{equation}
Combining the comparison inequality of Lemma~\ref{le:CosLawF} with
\eqref{eq:direct-inner-product-estimate} gives
\begin{align}
d^2(p_{k+1},p^*)
&\leq d^2(p_k,p^*)-d^2(p_{k+1},p_k)
+2\bigl\langle\log_{p_{k+1}}p_k,\log_{p_{k+1}}p^*\bigr\rangle
\notag\\
&\leq d^2(p_k,p^*)-d^2(p_{k+1},p_k)
+2\sigma_k\,d(p_{k+1},p_k)\,d(p_{k+1},p^*)
\label{eq:direct-multiplicative-fejer}\\
&=d^2(p_k,p^*)-\bigl(d(p_{k+1},p_k)-\sigma_k\,d(p_{k+1},p^*)\bigr)^2
+\sigma_k^2\,d^2(p_{k+1},p^*)
\notag\\
&\leq d^2(p_k,p^*)+\sigma_k^2\,d^2(p_{k+1},p^*).
\notag
\end{align}
Because \(0\leq\sigma_k\leq\bar\sigma<1\), this can be rewritten as
\begin{equation}
\label{eq:direct-quasi-fejer-multiplicative}
d^2(p_{k+1},p^*)
\leq\frac{1}{1-\sigma_k^2}\,d^2(p_k,p^*)
=\Bigl(1+\frac{\sigma_k^2}{1-\sigma_k^2}\Bigr)d^2(p_k,p^*)
\leq\Bigl(1+\frac{\sigma_k^2}{1-\bar\sigma^2}\Bigr)d^2(p_k,p^*).
\end{equation}
Iterating~\eqref{eq:direct-quasi-fejer-multiplicative} from \(k=0\) to \(N-1\)
and using \(1+t\leq e^{t}\) for \(t\geq0\) together with
Remark~\ref{rem:sigma-sq-summable}, we obtain, for every \(N\geq1\),
\[
d^2(p_N,p^*)
\leq d^2(p_0,p^*)\prod_{k=0}^{N-1}\Bigl(1+\frac{\sigma_k^2}{1-\bar\sigma^2}\Bigr)
\leq d^2(p_0,p^*)\,e^{\frac{1}{1-\bar\sigma^2}\sum_{k=0}^{N-1}\sigma_k^2}
\leq d^2(p_0,p^*)\,e^{\frac{\Sigma_2}{1-\bar\sigma^2}} .
\]
Taking square roots establishes (i). Note that the case \(N=0\) is trivial.

\medskip
\noindent(ii) By (i), the bound in (i) also holds with \(k+1\) in place of
\(k\), so the triangle inequality gives
\begin{equation}
\label{eq:direct-increment-bound}
d(p_{k+1},p_k)\leq d(p_{k+1},p^*)+d(p_k,p^*)
\leq2\,d(p_0,p^*)\,e^{\frac{\Sigma_2}{2(1-\bar\sigma^2)}} .
\end{equation}
On the other hand, rearranging~\eqref{eq:direct-multiplicative-fejer} yields
\[
d^2(p_{k+1},p_k)\leq d^2(p_k,p^*)-d^2(p_{k+1},p^*)
+2\sigma_k\,d(p_{k+1},p_k)\,d(p_{k+1},p^*),
\]
and, by \eqref{eq:direct-increment-bound} and (i),
\[
2\sigma_k\,d(p_{k+1},p_k)\,d(p_{k+1},p^*)
\leq2\sigma_k\cdot2\,d(p_0,p^*)\,e^{\frac{\Sigma_2}{2(1-\bar\sigma^2)}}
\cdot d(p_0,p^*)\,e^{\frac{\Sigma_2}{2(1-\bar\sigma^2)}}
=4\sigma_k\,d^2(p_0,p^*)\,e^{\frac{\Sigma_2}{1-\bar\sigma^2}} .
\]
Combining the last two displays proves (ii).

\medskip
\noindent(iii) By \eqref{eq:direct-proximal-inclusion},
\(v_k=\mu_k\log_{p_{k+1}}p_k-\varepsilon_k\); hence
\[
\|v_k\|\leq\mu_k\,\bigl\|\log_{p_{k+1}}p_k\bigr\|+\|\varepsilon_k\|
=\mu_k\,d(p_{k+1},p_k)+\|\varepsilon_k\| .
\]
Applying \eqref{eq:rockafellar-relative-error}, i.e.,
\(\|\varepsilon_k\|\leq\sigma_k\mu_k\,d(p_{k+1},p_k)\), and then
Assumption~\ref{ass:sigma-summable}, we get
$$
\|v_k\|\leq(1+\sigma_k)\mu_k\,d(p_{k+1},p_k)\leq(1+\bar\sigma)\mu_k\,d(p_{k+1},p_k).
$$
\end{proof}

Item~(i) of Lemma~\ref{le:direct-fejer-bounds} shows that the iterates remain
bounded, while item~(ii) shows that the increments \(d^2(p_{k+1},p_k)\) are
summable. Combined with the divergence of \(\sum_k\mu_k^{-2}\), this yields
convergence of the whole sequence to a point of \(S^*\).

\begin{theorem}
\label{thm:direct-inexact-proximal-convergence}
Let \(\{p_k\}_{k\in\mathbb N}\) be generated by Algorithm~\ref{Alg:SIPPM}.
Suppose that Assumption~\ref{ass:sigma-summable} holds and that
\(\sum_{k=0}^{+\infty}\mu_k^{-2}=+\infty\). Then
\begin{itemize}
\item[(i)] \(\displaystyle\sum_{k=0}^{+\infty}d^2(p_{k+1},p_k)<+\infty\);
\item[(ii)] there exists \(p^*\in S^*\) such that
\(\displaystyle\lim_{k\to+\infty}p_k=p^*\).
\end{itemize}
\end{theorem}

\begin{proof}
(i) Fix \(p^*\in S^*\). Summing Lemma~\ref{le:direct-fejer-bounds}(ii) from
\(k=0\) to \(N-1\) and using Assumption~\ref{ass:sigma-summable}, we obtain, for
every \(N\geq1\),
\[
\sum_{k=0}^{N-1}d^2(p_{k+1},p_k)
\leq d^2(p_0,p^*)-d^2(p_N,p^*)+4\Sigma\,d^2(p_0,p^*)\,e^{\frac{\Sigma_2}{1-\bar\sigma^2}}
\leq\Bigl(1+4\Sigma\,e^{\frac{\Sigma_2}{1-\bar\sigma^2}}\Bigr)d^2(p_0,p^*).
\]
The partial sums are nondecreasing and bounded above. Letting \(N\to+\infty\)
proves~(i).

\medskip
\noindent(ii) By~(i) and Lemma~\ref{le:direct-fejer-bounds}(iii),
\((\|v_k\|/\mu_k)^2\leq(1+\bar\sigma)^2d^2(p_{k+1},p_k)\), and therefore
\begin{equation}
\label{eq:direct-weighted-subgradient}
\sum_{k=0}^{+\infty}\frac{\|v_k\|^2}{\mu_k^2}<+\infty .
\end{equation}
We claim that \(\liminf_{k\to+\infty}\|v_k\|=0\). Otherwise, there would exist
\(\delta>0\) and \(k_0\in\mathbb N\) with \(\|v_k\|\geq\delta\) for all
\(k\geq k_0\), and the divergence of \(\sum_k\mu_k^{-2}\) would give
\[
\sum_{k=k_0}^{+\infty}\frac{\|v_k\|^2}{\mu_k^2}
\geq\delta^2\sum_{k=k_0}^{+\infty}\mu_k^{-2}=+\infty,
\]
contradicting~\eqref{eq:direct-weighted-subgradient}. Hence there is a
subsequence \(\{v_{k_\ell}\}_{\ell\in\mathbb N}\) with \(v_{k_\ell}\to0\). Set
\(m_\ell:=k_\ell+1\), so that \(v_{m_\ell-1}\in\partial f(p_{m_\ell})\) and
\(v_{m_\ell-1}\to0\).

By Lemma~\ref{le:direct-fejer-bounds}(i), \(\{p_k\}_{k\in\mathbb N}\) is
bounded. Since \(\mathbb M\) is complete and finite-dimensional, the
Hopf--Rinow theorem implies that bounded closed sets are compact. Thus, passing
to a further subsequence (not relabeled), we may assume \(p_{m_\ell}\to\tilde p\)
for some \(\tilde p\in\mathbb M\). Because
\(v_{m_\ell-1}\in\partial f(p_{m_\ell})\) and \(v_{m_\ell-1}\to0\),
Proposition~\ref{prop:subdifferential-continuity} yields
\(0\in\partial f(\tilde p)\), i.e., \(\tilde p\in S^*\).

It remains to show that the whole sequence converges to \(\tilde p\). Taking
\(p^*=\tilde p\) in~\eqref{eq:direct-quasi-fejer-multiplicative},
\[
d^2(p_{k+1},\tilde p)\leq\Bigl(1+\frac{\sigma_k^2}{1-\bar\sigma^2}\Bigr)d^2(p_k,\tilde p)
\;\;\text{for all }k\in\mathbb N .
\]
Iterating from \(k=m_{\ell_0}\) to \(m_{\ell_0}+N-1\), for any \(N\geq1\) and
\(\ell_0\in\mathbb N\), and using \(1+t\leq e^t\) and
Remark~\ref{rem:sigma-sq-summable}, we get
\[
d^2(p_{m_{\ell_0}+N},\tilde p)
\leq d^2(p_{m_{\ell_0}},\tilde p)\,
e^{\frac{1}{1-\bar\sigma^2}\sum_{k=m_{\ell_0}}^{+\infty}\sigma_k^2}.
\]
Let \(\varepsilon>0\). Since \(p_{m_\ell}\to\tilde p\) and \(\Sigma_2<+\infty\),
we may choose \(\ell_0\) so large that \(d^2(p_{m_{\ell_0}},\tilde p)<\varepsilon\)
and the exponential factor above is at most \(2\). Then
\(d^2(p_{m_{\ell_0}+N},\tilde p)<2\varepsilon\) for all \(N\geq1\), which proves
that \(d(p_k,\tilde p)\to0\). Setting \(p^*:=\tilde p\) completes the proof.
\end{proof}

\subsection{Iteration Complexity Bounds for Algorithm~\ref{Alg:SIPPM}}
\label{subsec:iter.comp}

We now derive iteration-complexity bounds for Algorithm~\ref{Alg:SIPPM}, both
for a stationarity measure and for function values. For a fixed integer
\(N\geq1\), let
\begin{equation}
\label{eq:mu-max}
\mu_N^{\max}:=\max_{0\leq k\leq N-1}\mu_k .
\end{equation}
We first bound the smallest
subgradient norm produced within the first \(N\) iterations.

\begin{theorem}
\label{thm:direct-complexity}
Let \(p^*\in S^*\), let \(\{p_k\}_{k\in\mathbb N}\) be an infinite sequence
generated by Algorithm~\ref{Alg:SIPPM}, and suppose that
Assumption~\ref{ass:sigma-summable} holds. Then, for every integer \(N\geq1\),
\[
\min_{0\leq k\leq N-1}\|v_k\|
\leq(1+\bar\sigma)\,\mu_N^{\max}\,d(p_0,p^*)
\sqrt{\frac{1+4\Sigma\,e^{\frac{\Sigma_2}{1-\bar\sigma^2}}}{N}} .
\]
\end{theorem}

\begin{proof}
By Lemma~\ref{le:direct-fejer-bounds}(iii), for \(0\leq k\leq N-1\),
\[
\|v_k\|^2\leq(1+\bar\sigma)^2\mu_k^2\,d^2(p_{k+1},p_k)
\leq(1+\bar\sigma)^2(\mu_N^{\max})^2\,d^2(p_{k+1},p_k).
\]
Summing over \(k=0,\dots,N-1\) gives
\[
\sum_{k=0}^{N-1}\|v_k\|^2
\leq(1+\bar\sigma)^2(\mu_N^{\max})^2\sum_{k=0}^{N-1}d^2(p_{k+1},p_k).
\]
Summing Lemma~\ref{le:direct-fejer-bounds}(ii) over the same range yields
\[
\sum_{k=0}^{N-1}d^2(p_{k+1},p_k)
\leq d^2(p_0,p^*)+4\,d^2(p_0,p^*)\,e^{\frac{\Sigma_2}{1-\bar\sigma^2}}\sum_{k=0}^{N-1}\sigma_k
\leq\Bigl(1+4\Sigma\,e^{\frac{\Sigma_2}{1-\bar\sigma^2}}\Bigr)d^2(p_0,p^*).
\]
Combining the two estimates,
\[
N\min_{0\leq k\leq N-1}\|v_k\|^2\leq\sum_{k=0}^{N-1}\|v_k\|^2
\leq(1+\bar\sigma)^2(\mu_N^{\max})^2\Bigl(1+4\Sigma\,e^{\frac{\Sigma_2}{1-\bar\sigma^2}}\Bigr)d^2(p_0,p^*).
\]
Dividing by \(N\) and taking square roots completes the proof.
\end{proof}

\begin{corollary}
\label{cor:direct-epsilon-complexity}
Under the assumptions of Theorem~\ref{thm:direct-complexity}, suppose that
\(0<\mu_k\leq\mu<+\infty\) for all \(k\in\mathbb N\). Then, for every
\(\epsilon>0\), Algorithm~\ref{Alg:SIPPM} generates a point \(p_{k+1}\) with
\(\|v_k\|\leq\epsilon\) for some \(0\leq k\leq N-1\) whenever
\[
N\geq\left\lceil\frac{(1+\bar\sigma)^2\mu^2\,d^2(p_0,p^*)
\bigl(1+4\Sigma\,e^{\frac{\Sigma_2}{1-\bar\sigma^2}}\bigr)}{\epsilon^2}\right\rceil .
\]
Thus, the worst-case stationarity complexity is \(\mathcal O(\epsilon^{-2})\).
\end{corollary}

We now turn to complexity in terms of objective values. Let \(p^*\in S^*\) be a
solution of~\eqref{eq:problem}. Since Algorithm~\ref{Alg:SIPPM} relies on
inexact subproblem solves, monotone decrease of \(\{f(p_k)\}_{k\in\mathbb N}\)
cannot be guaranteed in general; global functional rates are therefore stated
for the best iterate generated up to iteration \(N\). For any
\(v_k\in\partial f(p_{k+1})\), the subgradient inequality and the
Cauchy--Schwarz inequality give
\begin{equation}
\label{eq:direct-functional-residual}
f(p_{k+1})-f(p^*)
\leq-\langle v_k,\log_{p_{k+1}}p^*\rangle
\leq\|v_k\|\,d(p_{k+1},p^*).
\end{equation}
Thus, controlling the functional error requires simultaneous control of
\(\|v_k\|\) and of the distance \(d(p_{k+1},p^*)\). For the latter, we
introduce the following local growth condition.

\begin{assumption}
\label{ass:local-growth-condition}
Let \(p^*\in S^*\). There exist constants \(a>0\) and \(\tau>0\) such that \(d(p,p^*)\leq a\,\|v\|\) for all \(p\in\mathbb M\) and \(v\in\partial f(p)\) satisfying \(\|v\|\leq\tau\).
\end{assumption}

\begin{remark}
\label{rem:local-growth-condition}
Assumption~\ref{ass:local-growth-condition} implies that the minimizer is unique, that is, \(S^*=\{p^*\}\): indeed, if \(\tilde p\in S^*\), then \(0\in\partial f(\tilde p)\), which implies \(d(\tilde p,p^*)\leq a\|0\|=0\). Geometrically, this assumption is the Riemannian analogue of the Lipschitz-type behavior of \((\partial f)^{-1}\) at the origin studied by Rockafellar \cite[eq.~(3.1)]{Rockafellar76} and Solodov and Svaiter \cite[Theorem~2.2]{SolodovBenar1999}. A standard sufficient condition for Assumption~\ref{ass:local-growth-condition} is the \(m\)-strong convexity of \(f\). To see this, let \(p\in\mathbb M\) and \(v\in\partial f(p)\). Since \(0\in\partial f(p^*)\), inequality~\eqref{eq:strong-subgradient-inequality} yields
\[
f(p^*)\geq f(p)+\langle v,\log_p p^*\rangle+\frac m2\,d^2(p,p^*)
\quad\text{and}\quad
f(p)\geq f(p^*)+\frac m2\,d^2(p,p^*).
\]
Adding these two inequalities gives \(0\geq\langle v,\log_p p^*\rangle+m\,d^2(p,p^*)\). Applying the Cauchy--Schwarz inequality then leads to \(m\,d^2(p,p^*)\leq-\langle v,\log_p p^*\rangle\leq\|v\|\,d(p,p^*)\). Consequently, Assumption~\ref{ass:local-growth-condition} holds with \(a=1/m\) for any \(\tau>0\).
\end{remark}

\begin{theorem}
\label{thm:direct-functional-complexity}
Let \(p^*\in S^*\), suppose that Assumptions~\ref{ass:sigma-summable}
and~\ref{ass:local-growth-condition} hold, and let \(\{p_k\}_{k\in\mathbb N}\)
be generated by Algorithm~\ref{Alg:SIPPM}. If \(N\geq1\) satisfies
\begin{equation}
\label{eq:direct-N-large-enough}
(1+\bar\sigma)\,\mu_N^{\max}\,d(p_0,p^*)
\sqrt{\frac{1+4\Sigma\,e^{\frac{\Sigma_2}{1-\bar\sigma^2}}}{N}}\leq\tau,
\end{equation}
then
\[
\min_{0\leq k\leq N-1}\bigl[f(p_{k+1})-f(p^*)\bigr]
\leq\frac{a(1+\bar\sigma)^2(\mu_N^{\max})^2\,d^2(p_0,p^*)
\bigl(1+4\Sigma\,e^{\frac{\Sigma_2}{1-\bar\sigma^2}}\bigr)}{N}.
\]
\end{theorem}

\begin{proof}
Pick \(k_N\in\{0,\dots,N-1\}\) with
\(\|v_{k_N}\|=\min_{0\leq k\leq N-1}\|v_k\|\). By
Theorem~\ref{thm:direct-complexity} and~\eqref{eq:direct-N-large-enough},
\(\|v_{k_N}\|\leq\tau\), so Assumption~\ref{ass:local-growth-condition} applied
to \((p_{k_N+1},v_{k_N})\) gives \(d(p_{k_N+1},p^*)\leq a\|v_{k_N}\|\).
Inserting this into~\eqref{eq:direct-functional-residual} with \(k=k_N\) and
using Theorem~\ref{thm:direct-complexity}, we obtain
\[
f(p_{k_N+1})-f(p^*)\leq a\,\|v_{k_N}\|^2
\leq\frac{a(1+\bar\sigma)^2(\mu_N^{\max})^2\,d^2(p_0,p^*)
\bigl(1+4\Sigma\,e^{\frac{\Sigma_2}{1-\bar\sigma^2}}\bigr)}{N}.
\]
The proof is complete since
\(\min_{0\leq k\leq N-1}[f(p_{k+1})-f(p^*)]\leq f(p_{k_N+1})-f(p^*)\).
\end{proof}

\subsection{Convergence Rate for Algorithm~\ref{Alg:SIPPM}}
\label{sec:direct-rate}

In this section, we analyze the local rate of convergence of Algorithm~\ref{Alg:SIPPM} under a local growth condition on the objective function. The core strategy consists in decoupling the effect of the inexact subproblem resolution from the contraction properties of the underlying exact proximal point mapping. 
To this end, for each \(k\in\mathbb N\), let \(p^{\mathrm{ex}}_{k+1}\) denote the  exact proximal point of \(p_k\) defined by~\eqref{eq:prox-exact}. By the optimality conditions of the proximal subproblem, we immediately obtain the characterization
\begin{equation}
\label{eq:exact-proximal-subgradient}
\mu_k\log_{p^{\mathrm{ex}}_{k+1}}p_k   \in\partial f(p^{\mathrm{ex}}_{k+1}).
\end{equation}
The following key lemma serves as a bridge connecting the actual, inexactly generated point \(q_k\) with its exact theoretical counterpart \(p^{\mathrm{ex}}_{k+1}\).

\begin{lemma}
\label{lem:exact-vs-inexact}
Let \((q_k,v_k,\varepsilon_k)\) satisfy~\eqref{eq:direct-proximal-inclusion}.
Then
\begin{equation}
\label{eq:inexact-exact-proximal-comparison}
d(q_k,p^{\mathrm{ex}}_{k+1})\leq\frac{\|\varepsilon_k\|}{\mu_k}.
\end{equation}
Moreover, for every \(p^*\in S^*\),
\begin{equation}
\label{eq:exact-proximal-pythagorean}
d^2(p_k,p^*)\geq d^2(p_k,p^{\mathrm{ex}}_{k+1})+d^2(p^{\mathrm{ex}}_{k+1},p^*).
\end{equation}
\end{lemma}

\begin{proof}
Fix \(x\in\mathbb M\). Since \(v_k\in\partial f(q_k)\),
\(f(x)\geq f(q_k)+\langle v_k,\log_{q_k}x\rangle\). Moreover, Lemma~\ref{le:CosLawF}
gives
\(d^2(x,p_k)\geq d^2(q_k,x)+d^2(q_k,p_k)-2\langle\log_{q_k}x,\log_{q_k}p_k\rangle\).
Multiplying the latter by \(\mu_k/2>0\) and adding it to the former, we get
\[
f(x)+\frac{\mu_k}{2}d^2(x,p_k)
\geq f(q_k)+\frac{\mu_k}{2}d^2(q_k,p_k)
+\bigl\langle v_k-\mu_k\log_{q_k}p_k,\ \log_{q_k}x\bigr\rangle
+\frac{\mu_k}{2}d^2(q_k,x).
\]
Define \(\phi(x):=f(x)+(\mu_k/2)d^2(x,p_k)\). By
\eqref{eq:direct-proximal-inclusion}, \(v_k-\mu_k\log_{q_k}p_k=-\varepsilon_k\),
so
\begin{equation}
\label{eq:phi-lower-bound-qk}
\phi(x)\geq\phi(q_k)-\langle\varepsilon_k,\log_{q_k}x\rangle+\frac{\mu_k}{2}d^2(q_k,x).
\end{equation}
Since \(p^{\mathrm{ex}}_{k+1}\) minimizes the \(\mu_k\)-strongly convex function \(\phi\),
Lemma~\ref{lem:cat0}, specifically~\eqref{eq:strong-min}, yields
\begin{equation}
\label{eq:phi-lower-bound-sk}
\phi(x)\geq\phi(p^{\mathrm{ex}}_{k+1})+\frac{\mu_k}{2}d^2(x,p^{\mathrm{ex}}_{k+1}).
\end{equation}
Setting \(x=p^{\mathrm{ex}}_{k+1}\) in~\eqref{eq:phi-lower-bound-qk} and \(x=q_k\)
in~\eqref{eq:phi-lower-bound-sk} gives
\begin{align*}
\phi(p^{\mathrm{ex}}_{k+1})&\geq\phi(q_k)-\langle\varepsilon_k,\log_{q_k}p^{\mathrm{ex}}_{k+1}\rangle+\frac{\mu_k}{2}d^2(q_k,p^{\mathrm{ex}}_{k+1}),\\
\phi(q_k)&\geq\phi(p^{\mathrm{ex}}_{k+1})+\frac{\mu_k}{2}d^2(q_k,p^{\mathrm{ex}}_{k+1}).
\end{align*}
Adding these inequalities and using the Cauchy--Schwarz inequality,
\[
\mu_k\,d^2(q_k,p^{\mathrm{ex}}_{k+1})\leq\langle\varepsilon_k,\log_{q_k}p^{\mathrm{ex}}_{k+1}\rangle
\leq\|\varepsilon_k\|\,d(q_k,p^{\mathrm{ex}}_{k+1}).
\]
If \(q_k=p^{\mathrm{ex}}_{k+1}\), \eqref{eq:inexact-exact-proximal-comparison} is trivial.
Otherwise, divide by \(\mu_k\,d(q_k,p^{\mathrm{ex}}_{k+1})>0\).

To prove~\eqref{eq:exact-proximal-pythagorean}, note that \(p^*\) minimizes \(f\). Using~\eqref{eq:exact-proximal-subgradient} and the subgradient inequality at \(p^*\), we obtain
\(\mu_k\langle\log_{p^{\mathrm{ex}}_{k+1}}p_k,\log_{p^{\mathrm{ex}}_{k+1}}p^*\rangle\leq f(p^*)-f(p^{\mathrm{ex}}_{k+1})\leq0\).
Together with Lemma~\ref{le:CosLawF},
\[
d^2(p^{\mathrm{ex}}_{k+1},p_k)+d^2(p^{\mathrm{ex}}_{k+1},p^*)
\leq d^2(p^{\mathrm{ex}}_{k+1},p_k)+d^2(p^{\mathrm{ex}}_{k+1},p^*)-2\langle\log_{p^{\mathrm{ex}}_{k+1}}p_k,\log_{p^{\mathrm{ex}}_{k+1}}p^*\rangle
\leq d^2(p_k,p^*),
\]
which completes the proof.
\end{proof}

\begin{theorem}
\label{thm:direct-local-linear-rate}
Let \(\{p_k\}_{k\in\mathbb N}\) be an infinite sequence generated by
Algorithm~\ref{Alg:SIPPM}. Suppose that Assumptions~\ref{ass:sigma-summable}
and~\ref{ass:local-growth-condition} hold and that
\(0<\mu_k\leq\bar\mu<+\infty\) for all \(k\). Then \(\{p_k\}\) converges to
\(p^*\), and there exists \(k_0\in\mathbb N\) such that
\begin{equation}
\label{eq:direct-local-linear-estimate}
d(p_{k+1},p^*)\leq\theta_k\,d(p_k,p^*)\quad\forall\,k\geq k_0,
\;\;\text{where}\quad
\theta_k:=\frac{a\mu_k}{\sqrt{1+a^2\mu_k^2}}+\frac{\sigma_k}{1-\sigma_k}.
\end{equation}
Moreover,
\begin{equation}
\label{eq:direct-uniform-factor1}
\limsup_{k\to\infty}\theta_k\leq\widehat\theta:=\frac{a\bar\mu}{\sqrt{1+a^2\bar\mu^2}}<1 .
\end{equation}
In particular, \(\{p_k\}\) converges \(Q\)-linearly to \(p^*\). If, in addition,
\(\mu_k\to0\), then \(\theta_k\to0\) and the convergence is \(Q\)-superlinear.
\end{theorem}
\begin{proof}
Since \(\mu_k\leq\bar\mu\), we have \(\mu_k^{-2}\geq\bar\mu^{-2}>0\), which implies that \(\sum_k\mu_k^{-2}=+\infty\). Theorem~\ref{thm:direct-inexact-proximal-convergence} then ensures that \(p_k\to\tilde p\in S^*\), while Assumption~\ref{ass:local-growth-condition} guarantees the uniqueness of the minimizer, i.e., \(S^*=\{p^*\}\) (see Remark~\ref{rem:local-growth-condition}). Thus, \(p_k\to p^*\). Consequently, there exists an index \(k_0\in\mathbb N\) such that 
\begin{equation}
    \bar\mu\,d(p_k,p^*)\leq\tau, \;\; \forall k\geq k_0.
    \label{eq:distance-bound}
\end{equation}
We first show that the infinite sequence satisfies \(p_k\neq p^*\) for all \(k\). Indeed, suppose \(p_k=p^*\) for some \(k\). Then \(0\in\partial f(p_k)\), which implies that \(p_k\) minimizes \(f(\cdot)+(\mu_k/2)d^2(\cdot,p_k)\), i.e., \(p^{\mathrm{ex}}_{k+1}=p_k\). Combining~\eqref{eq:inexact-exact-proximal-comparison} with~\eqref{eq:rockafellar-relative-error} yields
$
d(q_k,p_k) = d(q_k,p^{\mathrm{ex}}_{k+1}) \leq \sigma_k\,d(q_k,p_k).
$
Since \(0\leq\sigma_k<1\), this forces \(q_k=p_k\), meaning that the algorithm would have terminated at iteration \(k\), which contradicts the assumption that \(\{p_k\}\) is an infinite sequence.

From \eqref{eq:exact-proximal-pythagorean} and \eqref{eq:distance-bound}, it follows that
$
\mu_k\,d(p^{\mathrm{ex}}_{k+1},p_k) \leq \bar\mu\,d(p_k,p^*) \leq \tau$ for all $k\geq k_0.
$
In view of \eqref{eq:exact-proximal-subgradient} and the identity \(\|\mu_k\log_{p^{\mathrm{ex}}_{k+1}}p_k\| = \mu_k\,d(p^{\mathrm{ex}}_{k+1},p_k)\), Assumption~\ref{ass:local-growth-condition} gives
$
d(p^{\mathrm{ex}}_{k+1},p^*) \leq a\mu_k\,d(p^{\mathrm{ex}}_{k+1},p_k).
$
Substituting this bound into \eqref{eq:exact-proximal-pythagorean} leads to
\[
\Bigl(1+\frac{1}{a^2\mu_k^2}\Bigr) d^2(p^{\mathrm{ex}}_{k+1},p^*) \leq d^2(p_k,p^{\mathrm{ex}}_{k+1}) + d^2(p^{\mathrm{ex}}_{k+1},p^*) \leq d^2(p_k,p^*),
\]
which proves the contraction inequality
\begin{equation}
\label{eq:exact-proximal-contraction}
d(p^{\mathrm{ex}}_{k+1},p^*) \leq \frac{a\mu_k}{\sqrt{1+a^2\mu_k^2}}\,d(p_k,p^*), \;\; \forall\,k\geq k_0.
\end{equation}

To account for the effect of inexact subproblem resolution, we combine \eqref{eq:inexact-exact-proximal-comparison} and \eqref{eq:rockafellar-relative-error} to deduce
\begin{equation}
\label{eq:qk-sk-bound}
d(q_k,p^{\mathrm{ex}}_{k+1}) \leq \frac{\|\varepsilon_k\|}{\mu_k} \leq \sigma_k\,d(q_k,p_k).
\end{equation}
Applying the triangle inequality together with the bound \(d(p^{\mathrm{ex}}_{k+1},p_k)\leq d(p_k,p^*)\) from \eqref{eq:exact-proximal-pythagorean}, it follows that
$
d(q_k,p_k) \leq d(q_k,p^{\mathrm{ex}}_{k+1}) + d(p^{\mathrm{ex}}_{k+1},p_k) \leq \sigma_k\,d(q_k,p_k) + d(p_k,p^*).
$
Since \(0\leq\sigma_k<1\), rearranging terms gives
\begin{equation}
\label{eq:qk-pk-bound}
d(q_k,p_k) \leq \frac{d(p_k,p^*)}{1-\sigma_k}.
\end{equation}
Recalling that \(p_{k+1}=q_k\), we apply the triangle inequality once more and combine \eqref{eq:qk-sk-bound}, \eqref{eq:exact-proximal-contraction}, and \eqref{eq:qk-pk-bound} to obtain, for all \(k\geq k_0\),
\[
d(p_{k+1},p^*) \leq d(q_k,p^{\mathrm{ex}}_{k+1}) + d(p^{\mathrm{ex}}_{k+1},p^*) \leq \sigma_k\,d(q_k,p_k) + \frac{a\mu_k}{\sqrt{1+a^2\mu_k^2}}\,d(p_k,p^*) \leq \theta_k\,d(p_k,p^*),
\]
which establishes \eqref{eq:direct-local-linear-estimate}.
Finally, observe that the map \(t\mapsto at/\sqrt{1+a^2t^2}\) is strictly increasing on \([0,\infty)\). Hence, \(\mu_k\leq\bar\mu\) implies that \(a\mu_k/\sqrt{1+a^2\mu_k^2}\leq\widehat\theta\) for all \(k\). Moreover, \(\sigma_k/(1-\sigma_k)\to0\) because \(\sigma_k\to0\) by Assumption~\ref{ass:sigma-summable}. This proves \eqref{eq:direct-uniform-factor1}, i.e., \(\limsup_k\theta_k\leq\widehat\theta<1\). Consequently, for any \(\theta\in(\widehat\theta,1)\), there exists an index \(k_1\geq k_0\) such that \(\theta_k\leq\theta\) for all \(k\geq k_1\). Since \(p_k\neq p^*\) for all \(k\), dividing \eqref{eq:direct-local-linear-estimate} by \(d(p_k,p^*)\) confirms \(Q\)-linear convergence. If, in addition, \(\mu_k\to0\), both terms comprising \(\theta_k\) vanish as \(k\to\infty\), which gives \(d(p_{k+1},p^*)/d(p_k,p^*)\leq\theta_k\to0\), thereby establishing \(Q\)-superlinear convergence.
\end{proof}

\begin{remark}
\label{rem:explicit-k0-direct}
The index \(k_0\) in Theorem~\ref{thm:direct-local-linear-rate} marks the first
stage from which the error bound of
Assumption~\ref{ass:local-growth-condition} can be applied to the exact
proximal pairs \((p^{\mathrm{ex}}_{k+1},u^{\mathrm{ex}}_{k+1})\), \(u^{\mathrm{ex}}_{k+1}\in\partial f(p^{\mathrm{ex}}_{k+1})\). As the proof shows, any
\(k_0\) such that \(\bar\mu\,d(p_k,p^*)\leq\tau\) for all \(k\geq k_0\) works,
since then \(\|u^{\mathrm{ex}}_{k+1}\|=\mu_k\,d(p^{\mathrm{ex}}_{k+1},p_k)\leq\bar\mu\,d(p_k,p^*)\leq\tau\). Such a
\(k_0\) exists because \(p_k\to p^*\), but it is not computable in general: the
direct method is only quasi-Fej\'er, and \(\{d(p_k,p^*)\}\) need not be
monotone.

Nevertheless, \(k_0\) can be made explicit \emph{a priori}. By
Lemma~\ref{le:direct-fejer-bounds}(i),
\[
d(p_k,p^*)\leq R_{\rm D}:=d(p_0,p^*)\,e^{\Sigma_2/(2(1-\bar\sigma^2))},
\;\;\forall\,k\in\mathbb N,
\]
so that, whenever
\begin{equation}
\label{eq:k0-zero-condition}
\bar\mu\,d(p_0,p^*)\,e^{\Sigma_2/(2(1-\bar\sigma^2))}\leq\tau,
\end{equation}
one may take \(k_0=0\), and \(d(p_{k+1},p^*)\leq\theta_k\,d(p_k,p^*)\) holds from
the first iteration. Condition~\eqref{eq:k0-zero-condition} is a smallness
requirement on \(d(p_0,p^*)\) relative to \(\tau/\bar\mu\). It holds
automatically when \(\tau=+\infty\), as is the case under strong convexity
(Remark~\ref{rem:local-growth-condition}). In that situation the linear
rate is global:
\[
d(p_k,p^*)\leq d(p_0,p^*)\prod_{j=0}^{k-1}\theta_j,\;\;\forall\,k\in\mathbb N .
\]
\end{remark}


\section{Inexact Hybrid Proximal Point Method} \label{sec:proxLine}
We now present a hybrid projection--proximal strategy adapted to Hadamard manifolds for solving \eqref{eq:problem}, drawing inspiration from the framework of Solodov and Svaiter \cite{SolodovBenar1999}. Unlike the direct scheme (Algorithm~\ref{Alg:SIPPM}), the approximate proximal point \(q_k\) is not accepted directly as the next iterate \(p_{k+1}\). As emphasized in Remark~\ref{rem:fixed-error-failure}, setting \(p_{k+1}=q_k\) under a non-vanishing relative error parameter \(\sigma>0\) may prevent convergence. To overcome this limitation, we construct a separating half-space \(\mathcal{S}_{q_k,v_k}\) associated with \(v_k \in \partial f(q_k)\) and project \(p_k\) onto it. This projection step guarantees a Fej\'er-monotone decrease of the distance to the solution set with a strict residual term, even when the proximal subproblem is solved with a fixed tolerance parameter \(\sigma \in [0,1)\).

The precise steps of the proposed approach are detailed in Algorithm~\ref{Alg:HHPPM}.

\begin{algorithm}[H]
\begin{footnotesize}
\begin{description}
  \item[Step 0.]
        Choose \(p_{0}\in{\mathbb M}_{}\) and one parameter  \(\sigma \in[0,1)\);
        set the iteration counter \({k}\gets 0\).
  \item[Step 1.]
        Take  \(\mu_{k}>0\) and compute a triple \(\bigl(q_{k},v_{k},\varepsilon_{k}\bigr)\) with \( v_{k},\varepsilon_{k} \in T_{q_k} {\mathbb M}\) satisfying the conditions
        \begin{equation} \label{eq:step1}
            0 = v_{k}-\mu_{k}\log_{q_{k}}p_{k}+\varepsilon_{k},\quad \quad  v_{k}\in \partial f(q_{k}), 
        \end{equation} 
 with  the relative error  \(\varepsilon_{k}\)  satisfying 
        \begin{equation} \label{eq:step1Error}   
        \|\varepsilon_{k}\|  \leq  \sigma \max\left\{\mu_k d(q_k,p_k),\|v_k\| \right\};
          \end{equation} 
        \textbf{Stopping test:} if \(v_{k}=0\) or \(q_{k}=p_{k}\), \textbf{stop} and return \(q_{k}\). 
\item[Step 2.]      
        Set the next iterate as
         \begin{equation} \label{eq:stepNextIt}
            p_{k+1}:={\cal P}_{\mathcal S_{q_k, v_k}}(p_{k})
           \end{equation}
           where
               \begin{equation} \label{eq:step2}
            \mathcal S_{q_k, v_k} :=
            \bigl\{p\in{\mathbb M} : \langle v_{k}  , \log_{q_{k}} p \rangle \leq 0\bigr\};
          \end{equation} 
  \item[Step 3.]
        Set \({k}\gets {k+1}\) and go to \textbf{Step~1}.
\end{description}
\caption{Inexact Hybrid Proximal Point Method (IHPPM)}
\label{Alg:HHPPM}
\end{footnotesize}
\end{algorithm}


\medskip

Algorithm~\ref{Alg:HHPPM} offers an alternative to the exact and inexact proximal schemes previously discussed. Rather than solving the proximal subproblem exactly (as in Algorithm~\ref{Alg:ExactPPM}) or accepting \(q_k\) directly as the next iterate (as in Algorithm~\ref{Alg:SIPPM}), it incorporates a projection step as a geometric correction. Consequently, the method allows a constant error tolerance \(\sigma\in[0,1)\), bypassing the summability requirement on the error parameters in Algorithm~\ref{Alg:SIPPM}.
Crucially, the relative error criterion \eqref{eq:step1Error} in Step~1 serves as the natural Riemannian generalization of \cite[Eq.~(1.7)]{SolodovBenar1999}. It keeps the subgradient norm \(\|v_k\|\) and the proximal displacement \(\mu_k d(q_k,p_k)\) on comparable scales. Indeed, by \eqref{eq:step1}, we have \(\mu_k\log_{q_k}p_k=v_k+\varepsilon_k\). Since \(\|\log_{q_k}p_k\|=d(q_k,p_k)\), the reverse triangle inequality yields
$
\bigl|\|v_k\| - \mu_k d(q_k,p_k)\bigr| = \bigl|\|v_k\| - \|v_k+\varepsilon_k\|\bigr| \leq \|\varepsilon_k\|.
$
Combining this estimate with \eqref{eq:step1Error} and using the identity \(\max\{a,b\}-\min\{a,b\}=|a-b|\) for \(a,b\geq 0\), we obtain
\begin{equation}
\min\left\{\mu_k d(q_k,p_k),\|v_k\|\right\} \geq (1-\sigma) \max\left\{\mu_k d(q_k,p_k),\|v_k\|\right\}.
\label{eq:min-max-comparison}
\end{equation}
Consequently, whenever \(\mu_k d(q_k,p_k)\neq 0\), it follows that \(1-\sigma \leq \|v_k\| / (\mu_k d(q_k,p_k)) \leq 1/(1-\sigma)\).
Moreover, \eqref{eq:min-max-comparison} shows that \(\|v_k\|=0\) if and only if \(d(q_k,p_k)=0\). Thus, condition \eqref{eq:step1Error} prevents the subgradient norm and the proximal displacement from having disparate scales, adapting dynamically to the magnitude of the current proximal inclusion while preserving the quantitative separation required in the subsequent projection step.

Notice that when \(\mathbb{M}\) is the hyperbolic space \(\mathbb{H}^n\), Lemma~\ref{lem:projection-hyperplane} in Appendix~\ref{sec:hyperbolic-projection} provides an explicit, closed-form formula for computing the metric projection in \eqref{eq:stepNextIt}.
To ensure that the projection step in \eqref{eq:stepNextIt} is well posed in general, the separating set \(\mathcal{S}_{q_k,v_k}\) must possess appropriate convexity properties. Note that \(\mathcal{S}_{q_k,v_k}\) always contains the entire solution set \(S^*\): for any \(p^* \in S^*\), the minimality of \(p^*\) and the geodesic convexity of \(f\) imply \(\left\langle v_k, \log_{q_k} p^* \right\rangle \leq f(p^*) - f(q_k) \leq 0\), establishing \(S^* \subseteq \mathcal{S}_{q_k,v_k}\). Accordingly, throughout our convergence analysis, we adopt the following standard assumption.

\begin{assumption}\label{ass:conv.su.}
For every \(k\in\mathbb{N}\), the set \(\mathcal{S}_{q_k,v_k}\) defined in \eqref{eq:step2} is convex.
\end{assumption}

By continuity of \(p \mapsto \langle v_k, \log_{q_k} p \rangle\), the set \(\mathcal{S}_{q_k,v_k}\) is closed. It is also non-empty since \(S^* \subseteq \mathcal{S}_{q_k,v_k}\). Under Assumption~\ref{ass:conv.su.}, \(\mathcal{S}_{q_k,v_k}\) is a non-empty, closed, convex subset of \(\mathbb{M}\), guaranteeing that the projection in Step~2 of Algorithm~\ref{Alg:HHPPM} is uniquely defined. Sufficient geometric conditions on the underlying manifold \(\mathbb{M}\) that guarantee the validity of Assumption~\ref{ass:conv.su.} are thoroughly discussed in Appendix~\ref{App:ass:conv.su.}.

We conclude this section by verifying that the stopping criterion in Step~1 triggers if and only if a solution to problem~\eqref{eq:problem} has been found.

\begin{lemma}\label{lem:stopping-test}
If Algorithm~\ref{Alg:HHPPM} terminates at Step~1, then the returned point is a solution of problem~\eqref{eq:problem}.
\end{lemma}

\begin{proof}
If \(v_k=0\), then \(0\in\partial f(q_k)\), which implies that \(q_k\) solves \eqref{eq:problem} by convexity. Now suppose \(q_k=p_k\). Then \eqref{eq:step1} reduces to \(v_k=-\varepsilon_k\), and \eqref{eq:step1Error} yields
$
\|v_k\| = \|\varepsilon_k\| \leq \sigma\|v_k\|.
$
Since \(0\leq\sigma<1\), this inequality forces \(\|v_k\|=0\), i.e., \(v_k=0\). Hence \(0\in\partial f(q_k)\), confirming that \(q_k=p_k\) solves problem~\eqref{eq:problem}.
\end{proof}

\subsection{Asymptotic Convergence Analysis for Algorithm~\ref{Alg:HHPPM}}
\label{subsec:asymptotic-convergence}

We now investigate the asymptotic behavior of the sequence \(\{p_k\}_{k\in\mathbb{N}}\) generated by Algorithm~\ref{Alg:HHPPM}. Throughout this subsection, we assume that the sequence is infinite and that the solution set \(S^*\) of problem~\eqref{eq:problem} is non-empty. By Lemma~\ref{lem:stopping-test}, this ensures that \(v_k\neq 0\) and \(q_k\neq p_k\) for every \(k\in\mathbb{N}\).

We first isolate essential angle and norm bounds resulting from the relative error criterion \eqref{eq:step1Error}. These estimates depend only on which term in the maximum is active at iteration \(k\) and will be invoked repeatedly in the convergence analysis throughout this section.

%
%

\begin{lemma}
\label{lem:angle-estimates}
Let $k\in\mathbb{N}$ correspond to a nonterminal iteration of
Algorithm~\ref{Alg:HHPPM}. Then the following statements hold.
\begin{itemize}
\item[(i)] If $\mu_k d(q_k,p_k)\geq\|v_k\|$, then
\begin{equation} \label{eq:first-case-angle} \langle v_k,\log_{q_k}p_k\rangle \;\geq\; (1-\sigma)\mu_k d^2(q_k,p_k), \end{equation} \begin{equation} \label{eq:first-case-norm} \|v_k\| \;\leq\; (1+\sigma)\mu_k d(q_k,p_k). 
\end{equation}

\item[(ii)] If $\|v_k\|\geq\mu_k d(q_k,p_k)$, then
\begin{equation}
\label{eq:second-case-angle}
\langle v_k,\log_{q_k}p_k\rangle
\geq
\frac{1-\sigma}{\mu_k}\|v_k\|^2,
\end{equation}
\begin{equation}
\label{eq:second-case-norm}
\mu_k d(q_k,p_k)
\leq
(1+\sigma)\|v_k\|.
\end{equation}
\end{itemize}
\end{lemma}

\begin{proof}
\noindent (i)  The condition
$\mu_k d(q_k,p_k)\geq\|v_k\|$, together with
\eqref{eq:step1Error}, yields
$
\|\varepsilon_k\|
\leq
\sigma\mu_k d(q_k,p_k).
$
Therefore, by \eqref{eq:step1} and the Cauchy--Schwarz inequality,
\[
\begin{aligned}
\langle v_k,\log_{q_k}p_k\rangle
&=
\mu_k d^2(q_k,p_k)
-
\langle\varepsilon_k,\log_{q_k}p_k\rangle\\
&\geq
\mu_k d^2(q_k,p_k)
-
\|\varepsilon_k\|d(q_k,p_k)
\geq
(1-\sigma)\mu_k d^2(q_k,p_k),
\end{aligned}
\]
which proves \eqref{eq:first-case-angle}. Moreover,
\begin{align*}
\|v_k\|
&\leq
\|v_k+\varepsilon_k\|+\|\varepsilon_k\|\\
&\leq
\mu_k d(q_k,p_k)
+
\sigma\mu_k d(q_k,p_k)
=
(1+\sigma)\mu_k d(q_k,p_k).
\end{align*}
This establishes
\eqref{eq:first-case-norm}. 

\noindent (ii)  The assumption
$\|v_k\|\geq\mu_k d(q_k,p_k)$ and \eqref{eq:step1Error} imply
$
\|\varepsilon_k\|
\leq
\sigma\|v_k\|.
$
Using \eqref{eq:step1}, we  obtain
\begin{align*}
\langle v_k,\log_{q_k}p_k\rangle
&=
\frac{1}{\mu_k}
\langle v_k,v_k+\varepsilon_k\rangle\\
&\geq
\frac{1}{\mu_k}
\left(
\|v_k\|^2-\|v_k\|\,\|\varepsilon_k\|
\right) \geq
\frac{1-\sigma}{\mu_k}\|v_k\|^2,
\end{align*}
which gives \eqref{eq:second-case-angle}. Finally, by \eqref{eq:step1},
$
\mu_k d(q_k,p_k)
=
\|v_k+\varepsilon_k\|
\leq
\|v_k\|+\|\varepsilon_k\|
\leq
(1+\sigma)\|v_k\|,
$
and hence \eqref{eq:second-case-norm} follows.
\end{proof}

Equipped with Lemma~\ref{lem:angle-estimates}, we now establish the foundational Fej\'er monotonicity properties of \(\{p_k\}_{k\in\mathbb{N}}\).

\begin{lemma}
\label{lema-fejer}
Let $\mathbb{M}$ be a Hadamard manifold, and let
$\{v_k\}_{k\in\mathbb{N}}$, $\{q_k\}_{k\in\mathbb{N}}$, and
$\{p_k\}_{k\in\mathbb{N}}$ denote the sequences generated by
Algorithm~\ref{Alg:HHPPM}. If Assumption~\ref{ass:conv.su.} holds,
then the following properties are satisfied:
\begin{itemize}
\item[(i)] $\langle v_k,\log_{q_k}p_k\rangle
\geq
\dfrac{1-\sigma}{1+\sigma}\,\|v_k\|\,d(q_k,p_k)$ for all $k\in\mathbb{N}$;
\item[(ii)] $d(p_k,p_{k+1}) \geq \dfrac{\langle v_k,\log_{q_k}p_k\rangle}{\|v_k\|}
\geq
\dfrac{1-\sigma}{1+\sigma}\,d(q_k,p_k)$ for all $k\in\mathbb{N}$;
\item[(iii)] $d^2(p_{k+1},p^*)\leq d^2(p_k, p^*)-d^2(p_k,p_{k+1})$
for all $k\in\mathbb{N}$ and all $p^*\in S^*$. In particular, the
sequence $\{p_k\}_{k\in\mathbb{N}}$ is Fejér convergent with respect to
$S^*$ and, consequently, bounded.
\end{itemize}
\end{lemma}

\begin{proof}
\noindent (i) Fix $k\in\mathbb{N}$. By Lemma~\ref{lem:angle-estimates},
exactly one of the two alternatives $\mu_k d(q_k,p_k)\geq\|v_k\|$ or
$\|v_k\|\geq\mu_k d(q_k,p_k)$ occurs. If the first holds, combining
\eqref{eq:first-case-angle} and \eqref{eq:first-case-norm} gives
\[
\langle v_k,\log_{q_k}p_k\rangle
\;\geq\;
(1-\sigma)\mu_k d^2(q_k,p_k)
\;\geq\;
\frac{1-\sigma}{1+\sigma}\,\|v_k\|\,d(q_k,p_k).
\]
If instead the second holds, combining \eqref{eq:second-case-angle}
and \eqref{eq:second-case-norm} gives
\[
\langle v_k,\log_{q_k}p_k\rangle
\;\geq\;
\frac{1-\sigma}{\mu_k}\|v_k\|^2
\;\geq\;
\frac{1-\sigma}{1+\sigma}\,\|v_k\|\,d(q_k,p_k).
\]
In either case the claimed inequality holds.

\noindent (ii) By item~(i),
\begin{equation}
\label{eq:angle-lower-bound}
\langle v_k,\log_{q_k}p_k\rangle
\;\geq\;
\frac{1-\sigma}{1+\sigma}\,\|v_k\|\,d(q_k,p_k)
\;>\;0,
\end{equation}
where the strict inequality follows since $\sigma\in [0,1)$, $v_k\neq
0$, and $q_k\neq p_k$. Hence $p_k\notin\mathcal S_{q_k,v_k}$; see
\eqref{eq:step2}. Since $p_{k+1}\in\mathcal S_{q_k,v_k}$, we have
$\langle v_k,\log_{q_k}p_{k+1}\rangle\leq0$. Therefore,
\begin{align*}
\langle v_k,\log_{q_k}p_k\rangle
&=
\langle v_k,\log_{q_k}p_k-\log_{q_k}p_{k+1}\rangle
+
\langle v_k,\log_{q_k}p_{k+1}\rangle\\
&\leq
\|v_k\|\,\|\log_{q_k}p_k-\log_{q_k}p_{k+1}\|\\
&\leq
\|v_k\|\,d(p_k,p_{k+1}),
\end{align*}
where the last inequality follows from Lemma~\ref{le:CosLawF}.
Combining this estimate with \eqref{eq:angle-lower-bound} completes
the proof of (ii).

\noindent (iii) Let $p^*\in S^*$. By \eqref{eq:step1} and
Definition~\ref{def:subgradient},
$
\langle v_k,\log_{q_k}p^*\rangle
\;\leq\;
f(p^*)-f(q_k)
\;\leq\;0,
$
so $p^*\in\mathcal S_{q_k,v_k}$. Since $\mathcal S_{q_k,v_k}$ is
convex by Assumption~\ref{ass:conv.su.}, Lemma~\ref{le:CosLawF},
\eqref{eq:stepNextIt}, and Proposition~\ref{prop:proj-ineq} yield
\begin{align*}
d^2(p_{k+1},p^*)
&\leq
d^2(p_k,p^*)-d^2(p_k,p_{k+1})
+
2\langle\log_{p_{k+1}}p_k,\log_{p_{k+1}}p^*\rangle\\
&\leq
d^2(p_k,p^*)-d^2(p_k,p_{k+1}).
\end{align*}
Thus $\{p_k\}_{k\in\mathbb{N}}$ is Fejér convergent with respect to the
solution set and, consequently, bounded (see
Definition~\ref{def:fejer} and Proposition~\ref{prop:fejer}).
\end{proof}

We are now ready to state and derive the main convergence result.

\begin{theorem}\label{hybrid-convergence}
Let 
$\{p_k\}_{k\in\mathbb{N}}$ be the sequence generated by
Algorithm~\ref{Alg:HHPPM}. Suppose that Assumption~\ref{ass:conv.su.}
holds and that
$
\sum_{k=0}^{+\infty}\mu_k^{-2}=+\infty.
$
Then the following statements hold:
\begin{itemize}
\item[(i)] $\displaystyle\sum_{k=0}^{+\infty}d^2(q_k,p_k)<+\infty$;
\item[(ii)] there exists $p^*\in S^*$ such that
$\displaystyle\lim_{k\to+\infty}p_k=p^*$.
\end{itemize}
\end{theorem}

\begin{proof}
\noindent (i) Let $\bar p\in S^*$. By the Fejér-type inequality
established in Lemma~\ref{lema-fejer}(iii),
\begin{equation}
\label{eq:fejer-monotone-strict}
d^2(p_{k+1},\bar p)
\;\leq\;
d^2(p_k,\bar p)-d^2(p_k,p_{k+1}),
\;\;\forall k\in\mathbb{N}.
\end{equation}
Summing from $k=0$ to $N-1$ gives
\[
\sum_{k=0}^{N-1}d^2(p_k,p_{k+1})
\;\leq\;
d^2(p_0,\bar p)-d^2(p_N,\bar p)
\;\leq\;
d^2(p_0,\bar p).
\]
Letting $N\to+\infty$, we conclude that
$\sum_{k=0}^{+\infty}d^2(p_k,p_{k+1})<+\infty$. Hence, the desired inequality follows directly from Lemma~\ref{lema-fejer}(ii).

\noindent (ii) Since \eqref{eq:min-max-comparison} yields
$(1-\sigma)\|v_k\|\leq\mu_k d(q_k,p_k)$, item~(i) gives
\begin{equation}
\label{eq:weighted-subgradient-summable}
\sum_{k=0}^{+\infty}\frac{\|v_k\|^2}{\mu_k^2}
\;\leq\;
\frac{1}{(1-\sigma)^2}
\sum_{k=0}^{+\infty}d^2(q_k,p_k)
\;<\;+\infty.
\end{equation}
Suppose, for contradiction, that
$\liminf_{k\to+\infty}\|v_k\|>0$. Then there exist $\delta>0$ and
$k_0\in\mathbb{N}$ such that $\|v_k\|\geq\delta$ for all $k\geq k_0$,
so that
\[
\sum_{k=k_0}^{\infty}\frac{\|v_k\|^2}{\mu_k^2}
\;\geq\;
\delta^2\sum_{k=k_0}^{\infty}\frac{1}{\mu_k^2}
\;=\;
+\infty,
\]
contradicting \eqref{eq:weighted-subgradient-summable}. Hence
$\liminf_{k\to\infty}\|v_k\|=0$, and there is a subsequence
$\{v_{k_\ell}\}_{\ell\in\mathbb{N}}$ with $v_{k_\ell}\to0$.
Since $\{p_k\}_{k\in\mathbb{N}}$ is bounded and $\mathbb{M}$ is a
finite-dimensional complete Riemannian manifold, the Hopf--Rinow
theorem ensures that $\{p_{k_\ell}\}_{\ell\in\mathbb{N}}$ admits a
convergent subsequence; without relabeling, let $p_{k_\ell}\to p^*$.
Then $d(p_{k_\ell},q_{k_\ell})\to0$, so $q_{k_\ell}\to p^*$. By
Proposition~\ref{prop:subdifferential-continuity}, $0\in\partial
f(p^*)$, and therefore $p^*\in S^*$.
Finally, applying \eqref{eq:fejer-monotone-strict} with $\bar p=p^*$
shows that $\{d(p_k,p^*)\}_{k\in\mathbb{N}}$ is nonincreasing. Since
$d(p_{k_\ell},p^*)\to0$, we conclude $d(p_k,p^*)\to0$, i.e.,
$p_k\to p^*$, which completes the proof.
\end{proof}

\subsection{Iteration Complexity Bounds for Algorithm~\ref{Alg:HHPPM}}
\label{subsec:iter.comp.proj}

In this section, we derive complexity bounds for Algorithm~\ref{Alg:HHPPM} using $\mu_N^{\max}$ defined in~\eqref{eq:mu-max}. We first refine the projection-displacement estimate, sharpening the lower bound on $d(p_k, p_{k+1})$ from Lemma~\ref{lema-fejer}(ii) by incorporating both conditions of the relative error criterion in Lemma~\ref{lem:angle-estimates}.

\begin{lemma}
\label{lem:hybrid-complexity-step}
Let $k \in \mathbb{N}$ denote the index of a nonterminal iteration of Algorithm~\ref{Alg:HHPPM}. Then,
\begin{equation}
\label{eq:hybrid-sharp}
d(p_k, p_{k+1}) \ge (1-\sigma) \max\left\{ d(q_k, p_k), \, \frac{\|v_k\|}{\mu_k} \right\}.
\end{equation}
\end{lemma}

\begin{proof}
First, assume that $\mu_k d(q_k, p_k) \ge \|v_k\|$. By \eqref{eq:first-case-angle}, we have
\[
\langle v_k, \log_{q_k} p_k \rangle \ge (1-\sigma)\mu_k d^2(q_k, p_k).
\]
Combining this inequality with Lemma~\ref{lema-fejer}(ii) yields
\begin{align*}
d(p_k, p_{k+1}) 
&\ge \frac{\langle v_k, \log_{q_k} p_k \rangle}{\|v_k\|} 
\ge (1-\sigma) \frac{\mu_k d^2(q_k, p_k)}{\|v_k\|} \\
&\ge (1-\sigma) d(q_k, p_k) 
= (1-\sigma) \max\left\{ d(q_k, p_k), \, \frac{\|v_k\|}{\mu_k} \right\},
\end{align*}
where the last equality follows directly from the initial assumption.
Conversely, suppose that $\|v_k\| \ge \mu_k d(q_k, p_k)$. In view of \eqref{eq:second-case-angle},
$
\langle v_k, \log_{q_k} p_k \rangle \ge (1-\sigma)\mu_k^{-1} \|v_k\|^2.
$
Applying Lemma~\ref{lema-fejer}(ii) once again, we deduce that
\[
d(p_k, p_{k+1}) \ge \frac{\langle v_k, \log_{q_k} p_k \rangle}{\|v_k\|} \ge (1-\sigma) \frac{\|v_k\|}{\mu_k} = (1-\sigma) \max\left\{ d(q_k, p_k), \, \frac{\|v_k\|}{\mu_k} \right\},
\]
which establishes \eqref{eq:hybrid-sharp} and completes the proof.
\end{proof}

We can now obtain the iteration-complexity estimate for the hybrid method.

\begin{theorem}
\label{thm:hybrid-complexity}
Let $p^* \in S^*$, and let $\{p_k\}_{k\in\mathbb{N}}$ be an infinite sequence generated by Algorithm~\ref{Alg:HHPPM}, and suppose that Assumption~\ref{ass:conv.su.} holds. Then, for every integer $N \ge 1$,
\[
\min_{0 \le k \le N-1} \|v_k\| \le \frac{\mu_N^{\max} d(p_0, p^*)}{(1 - \sigma)\sqrt{N}}.
\]
\end{theorem}

\begin{proof}
By Lemma~\ref{lem:hybrid-complexity-step}, we have $d(p_k, p_{k+1}) \ge (1 - \sigma)\mu_k^{-1} \|v_k\|$ for all $k \ge 0$. Rearranging this inequality and taking the square yields
\[
\|v_k\|^2 \le \frac{\mu_k^2}{(1 - \sigma)^2} d^2(p_k, p_{k+1}) \le \frac{(\mu_N^{\max})^2}{(1 - \sigma)^2} d^2(p_k, p_{k+1}).
\]
Summing this relation over $k = 0, \dots, N-1$ and applying Lemma~\ref{lema-fejer}(iii), we obtain
\[
\sum_{k=0}^{N-1} \|v_k\|^2 \le \frac{(\mu_N^{\max})^2}{(1 - \sigma)^2} \sum_{k=0}^{N-1} d^2(p_k, p_{k+1}) \le \frac{(\mu_N^{\max})^2 d^2(p_0, p^*)}{(1 - \sigma)^2}.
\]
Consequently, we have
\[
N \min_{0 \le k \le N-1} \|v_k\|^2 \le \sum_{k=0}^{N-1} \|v_k\|^2 \le \frac{(\mu_N^{\max})^2 d^2(p_0, p^*)}{(1 - \sigma)^2}.
\]
Dividing by $N$ and taking the square root on both sides completes the proof.
\end{proof}

\begin{corollary}
\label{cor:hybrid-epsilon-complexity}
Under the assumptions of Theorem~\ref{thm:hybrid-complexity}, suppose that $0 < \mu_k \le \mu < +\infty$ for all $k \in \mathbb{N}$. Then, Algorithm~\ref{Alg:HHPPM} generates an $\epsilon$-critical point (in the sense of Definition~\ref{def:epsilon-critical-point}) among $q_0, \dots, q_{N-1}$ whenever
\[
N \ge \left\lceil \frac{\mu^2 d^2(p_0, p^*)}{(1-\sigma)^2 \epsilon^2} \right\rceil.
\]
Thus, the worst-case stationarity complexity is $\mathcal{O}(\epsilon^{-2})$.
\end{corollary}

Next, we derive an iteration-complexity estimate in terms of the objective values. Because the sequence $\{q_k\}_{k\in\mathbb{N}}$ generated by Algorithm~\ref{Alg:HHPPM} is computed via inexact proximal steps, the sequence of objective values $\{f(q_k)\}_{k\in\mathbb{N}}$ is not necessarily nonincreasing. Consequently, a global complexity bound for the functional residual is naturally expressed in terms of the best iterate. To control the distance to $p^*$ via subgradient information, we adopt the local growth condition stated in Assumption~\ref{ass:local-growth-condition}.

\begin{theorem}
\label{thm:hybrid-functional-complexity}
Let $p^* \in S^*$ and suppose that Assumptions~\ref{ass:conv.su.} and~\ref{ass:local-growth-condition} hold. Let $\{p_k\}_{k\in\mathbb{N}}$ be the sequence generated by Algorithm~\ref{Alg:HHPPM}. If
\begin{equation}
\label{eq:N-large-enough}
\frac{\mu_N^{\max} d(p_0, p^*)}{(1-\sigma)\sqrt{N}} \le \tau,
\end{equation}
then
\begin{equation}
\label{eq:hybrid-best-functional}
\min_{0 \le k \le N-1} \bigl[ f(q_k) - f(p^*) \bigr] 
\le 
\frac{a(\mu_N^{\max})^2 d^2(p_0, p^*)}{(1-\sigma)^2 N}.
\end{equation}
\end{theorem}

\begin{proof}
Let $k_N \in \{0, \dots, N-1\}$ be an index satisfying $\|v_{k_N}\| = \min_{0 \le k \le N-1} \|v_k\|$. By Theorem~\ref{thm:hybrid-complexity}, it holds that
\begin{equation}
\label{eq:best-subgradient-squared}
\|v_{k_N}\|^2 \le \frac{(\mu_N^{\max})^2 d^2(p_0, p^*)}{(1-\sigma)^2 N}.
\end{equation}
Condition~\eqref{eq:N-large-enough} ensures that $\|v_{k_N}\| \le \tau$, allowing us to apply Assumption~\ref{ass:local-growth-condition} to $(q_{k_N}, v_{k_N})$, which yields $d(q_{k_N}, p^*) \le a \|v_{k_N}\|$. On the other hand, for $v_{k_N} \in \partial f(q_{k_N})$, the subdifferential inequality together with Cauchy--Schwarz gives
\begin{equation}
\label{eq:functional-residual}
f(q_{k_N}) - f(p^*) \le -\langle v_{k_N}, \log_{q_{k_N}} p^* \rangle \le \|v_{k_N}\| \, d(q_{k_N}, p^*).
\end{equation}
Combining $d(q_{k_N}, p^*) \le a \|v_{k_N}\|$ with \eqref{eq:functional-residual} and \eqref{eq:best-subgradient-squared}, we obtain
\[
f(q_{k_N}) - f(p^*) \le a \|v_{k_N}\|^2 \le \frac{a(\mu_N^{\max})^2 d^2(p_0, p^*)}{(1-\sigma)^2 N}.
\]
The result now follows immediately since $\displaystyle\min_{0 \le k \le N-1} [f(q_k) - f(p^*)] \le f(q_{k_N}) - f(p^*)$.
\end{proof}

\begin{remark}
\label{rem:N-large-enough}
Condition~\eqref{eq:N-large-enough} should not be interpreted as an additional structural assumption on the problem. Its goal is only to guarantee that, among the first $N$ iterations, there exists an index $k_N\in\{0,\ldots,N-1\}$ for which the local growth condition can be applied. Indeed, from the previous complexity estimate,
$$
\min_{0\leq k\leq N-1}\|v_k\|
\leq
\frac{\mu_N^{\max}d(p_0,\bar p)}
{(1-\sigma)\sqrt N}.
$$
Hence, if~\eqref{eq:N-large-enough} holds, one can choose $k_N$ such that
$$
\|v_{k_N}\|
=
\min_{0\leq k\leq N-1}\|v_k\|
\leq \tau.
$$
Therefore, the local growth condition 
$
d(q_{k_N},p^*)\leq a\|v_{k_N}\|
$ 
is valid at this iterate, which allows us to derive the functional residual estimate. 
Moreover, if the sequence $\{\mu_k\}$ is bounded above, say $\mu_k\leq\bar\mu$ for every $k$, then condition~\eqref{eq:N-large-enough} is automatically satisfied for all sufficiently large $N$. In fact, it is enough to require
$$
N\geq\frac{\bar\mu^2d^2(p_0, p^*)} {(1-\sigma)^2\tau^2}.
$$
Thus, condition~\eqref{eq:N-large-enough} simply identifies the iteration range in which the local growth assumption becomes applicable.
\end{remark}

\subsection{Convergence Rate for Algorithm~\ref{Alg:HHPPM}}

We now establish a sharper local linear convergence estimate for Algorithm~\ref{Alg:HHPPM}. The proof considers separately the two cases determined by the maximum in the relative error criterion \eqref{eq:step1Error}, allowing us to exploit the specific estimates available in each case and obtain a sharper contraction factor.

\begin{theorem}
\label{thm:hybrid-local-linear-rate}
Let $\{p_k\}_{k\in\mathbb{N}}$ be an infinite sequence generated by Algorithm~\ref{Alg:HHPPM}. Suppose that Assumptions~\ref{ass:conv.su.} and~\ref{ass:local-growth-condition} hold at a solution $p^*$, and that $0<\mu_k\leq\bar{\mu}<+\infty$ for all $k$. Then $\{p_k\}_{k\in\mathbb{N}}$ converges to $p^*$, and, for all sufficiently large $k$,
\begin{equation}
\label{eq:hybrid-local}
d(p_{k+1},p^*)\leq\rho_k\,d(p_k,p^*) \quad \text{where} \quad \rho_k
:=
\left[
1-
\left(
\frac{1-\sigma}{1+a\mu_k}
\right)^2
\right]^{1/2}.
\end{equation}
Moreover,
\begin{equation}
\label{eq:direct-uniform-factor}
\limsup_{k\to\infty}\rho_k\ \leq\ \widehat{\rho}:=\left[
1-
\left(
\frac{1-\sigma}{1+a\bar\mu}
\right)^2
\right]^{1/2}<1.
\end{equation}
In particular, $\{p_k\}$ converges $Q$-linearly to $p^*$. 
\end{theorem}
\begin{proof}
Since \(\mu_k \le \bar{\mu}\) for all \(k\), we have \(\sum_{k} \mu_k^{-2} \ge \sum_{k} \bar{\mu}^{-2} = +\infty\). Consequently, Theorem~\ref{hybrid-convergence} applies, yielding \(\sum_k d^2(q_k, p_k) < +\infty\), which in turn implies \(d(q_k, p_k) \to 0\) and \(p_k \to p^*\) for some \(p^* \in S^*\). Furthermore, Assumption~\ref{ass:local-growth-condition} ensures that \(S^*\) is a singleton, i.e., \(S^* = \{p^*\}\) (see Remark~\ref{rem:local-growth-condition}).

Next, we show that \(\|v_k\| \to 0\). If \(\|v_k\| \ge \mu_k d(q_k, p_k)\), then \(\|\varepsilon_k\| \le \sigma \|v_k\|\) by virtue of \eqref{eq:step1Error}. Combining this with the reverse triangle inequality
$
\mu_k d(q_k, p_k) = \|v_k + \varepsilon_k\| \ge \|v_k\| - \|\varepsilon_k\|,
$
we readily obtain
$
\|v_k\| \le \mu_k/(1 - \sigma) d(q_k, p_k).
$
Conversely, if \(\|v_k\| < \mu_k d(q_k, p_k)\), the same upper bound holds trivially since \(\sigma \in [0,1)\). Hence,
$
\|v_k\| \le \mu_k/(1 - \sigma) d(q_k, p_k) \to 0,
$
taking into account that \(\mu_k \le \bar{\mu}\) and \(d(q_k, p_k) \to 0\). Consequently, invoking Assumption~\ref{ass:local-growth-condition} for the pair \((q_k, v_k)\) for all sufficiently large \(k\) yields
$
d(q_k, p^*) \le a \|v_k\|.
$
By the triangle inequality,
\begin{equation}
\label{eq:hybrid-common-bound}
d(p_k,p^*)\leq d(p_k,q_k)+d(q_k,p^*)\leq d(q_k,p_k)+a\|v_k\|,
\end{equation}
a bound we now use in both cases of Lemma~\ref{lem:angle-estimates}.

\medskip
\noindent\textit{Case $\mu_kd(q_k,p_k)\geq\|v_k\|$.} By Lemma~\ref{lema-fejer}(ii) and \eqref{eq:first-case-angle} in Lemma~\ref{lem:angle-estimates}(i),
$$
d(p_k,p_{k+1})
\geq
\frac{\langle v_k,\log_{q_k}p_k\rangle}{\|v_k\|}
\geq
(1-\sigma)\frac{\mu_kd(q_k,p_k)}{\|v_k\|}d(q_k,p_k)
\geq
(1-\sigma)d(q_k,p_k),
$$
where the last step uses the case hypothesis $\mu_kd(q_k,p_k)\geq\|v_k\|$. The same hypothesis, substituted into \eqref{eq:hybrid-common-bound}, gives $d(p_k,p^*)\leq(1+a\mu_k)d(q_k,p_k)$, so that
\begin{equation}
\label{eq:hybrid-first-case}
d(p_k,p_{k+1})\geq\frac{1-\sigma}{1+a\mu_k}d(p_k,p^*).
\end{equation}

\medskip
\noindent\textit{Case $\|v_k\|\geq\mu_kd(q_k,p_k)$.} By Lemma~\ref{lema-fejer}(ii) and \eqref{eq:second-case-angle} in Lemma~\ref{lem:angle-estimates}(ii),
$$
d(p_k,p_{k+1})\geq\frac{\langle v_k,\log_{q_k}p_k\rangle}{\|v_k\|}\geq\frac{1-\sigma}{\mu_k}\|v_k\|.
$$
Setting $t_k:=\|v_k\|/(\mu_kd(q_k,p_k))\geq1$, this reads $d(p_k,p_{k+1})\geq(1-\sigma)t_kd(q_k,p_k)$. Substituting $\|v_k\|=t_k\mu_kd(q_k,p_k)$ into \eqref{eq:hybrid-common-bound} gives $d(p_k,p^*)\leq(1+a\mu_kt_k)d(q_k,p_k)$, so that
$$
d(p_k,p_{k+1})\geq\frac{(1-\sigma)t_k}{1+a\mu_kt_k}d(p_k,p^*).
$$
Since $t_k\geq1$, and this is equivalent to $t_k/(1+a\mu_kt_k)\geq1/(1+a\mu_k)$, we again obtain
\begin{equation}
\label{eq:hybrid-second-case}
d(p_k,p_{k+1})\geq\frac{1-\sigma}{1+a\mu_k}d(p_k,p^*).
\end{equation}

Combining \eqref{eq:hybrid-first-case} and \eqref{eq:hybrid-second-case}, we obtain, for all sufficiently large $k$,
\begin{equation}
\label{eq:hybrid-sharp-local}
d(p_k,p_{k+1})\geq\frac{1-\sigma}{1+a\mu_k}d(p_k,p^*).
\end{equation}
By Lemma~\ref{lema-fejer}(iii), together with \eqref{eq:hybrid-sharp-local}, we obtain
$$
d^2(p_{k+1},p^*)\leq d^2(p_k,p^*)-d^2(p_k,p_{k+1})\leq\left[1-\left(\frac{1-\sigma}{1+a\mu_k}\right)^2\right]d^2(p_k,p^*),
$$
and, after taking square roots, this proves \eqref{eq:hybrid-local}. 
Finally, since the function \(t \mapsto \left[ 1 - \left( \frac{1-\sigma}{1+at} \right)^2 \right]^{1/2}\) is strictly increasing on \([0, \infty)\) and \(\mu_k \le \bar{\mu}\), inequality \eqref{eq:direct-uniform-factor} follows directly, completing the proof.
\end{proof}

\begin{remark}
\label{rem:hybrid-superlinear}
A fixed relative tolerance \(\sigma \in (0,1)\) in Algorithm~\ref{Alg:HHPPM} guarantees local \(Q\)-linear convergence. However, it does not ensure \(Q\)-superlinear convergence even when \(\mu_k \to 0\), since the contraction factor satisfies
$
\lim_{\mu_k \to 0} \rho_k = \sqrt{1 - (1-\sigma)^2} > 0.
$
To achieve \(Q\)-superlinear convergence, Algorithm~\ref{Alg:HHPPM} can be implemented using adaptive tolerances \(\{\sigma_k\} \subset [0,1)\) such that \(\sigma_k \to 0\). In this setting, if \(\mu_k \to 0\), the corresponding contraction factor satisfies
\[
\rho_k = \left[ 1 - \left( \frac{1 - \sigma_k}{1 + a\mu_k} \right)^2 \right]^{1/2} \longrightarrow 0 \quad \text{as } k \to \infty,
\]
thereby establishing the \(Q\)-superlinear convergence of \(\{p_k\}\) to \(p^*\).
\end{remark}

\section{Conclusion}
\label{sec:conclusions}

In this paper, we investigated two inexact proximal point methods for solving convex optimization problems on Hadamard manifolds. 
The first approach is a standard inexact variant of the classical proximal point algorithm. It applies to general Hadamard manifolds under a standard summability condition on the relative errors. The second approach is a hybrid inexact proximal point method that incorporates a projection step onto a separating hyperplane-like set, allowing the summability requirement to be relaxed to a constant relative error tolerance. 
For both algorithms, we established full asymptotic convergence of the sequence of iterates and provided global iteration-complexity bounds to achieve approximate stationary points. Furthermore, under a local growth condition, we derived explicit complexity estimates for functional values along with local convergence rates. In particular, the direct method achieves a local \(Q\)-linear rate that accelerates to \(Q\)-superlinear when the proximal parameters tend to zero. For the hybrid method, we established local \(Q\)-linear convergence together with a refined contraction factor that explicitly incorporates the geometry of its dual error criterion.
From a geometric standpoint, we analyzed the geodesic convexity of the separating sets \(\mathcal{S}_{q_k,v_k}\) defined in \eqref{eq:step2}. We proved that these sets are automatically geodesically convex on Hadamard surfaces as well as on Hadamard manifolds of constant sectional curvature, thereby ensuring the structural validity of Assumption~\ref{ass:conv.su.} in those settings.
A promising direction for future research is to extend the geometric analysis of Assumption~\ref{ass:conv.su.} to broader classes of Hadamard manifolds with nonconstant sectional curvature, such as the manifold of Symmetric Positive Definite (SPD) matrices. A key avenue toward this goal is to investigate separation theorems \cite{silva2022fenchel} to establish precise geometric conditions under which the separating sets \(\mathcal{S}_{q_k,v_k}\) remain geodesically convex in nonconstant curvature environments.

\appendix

\section{Projection onto Hyperbolic Hyperplanes}
\label{sec:hyperbolic-projection}

In this section, we derive an explicit formula for the metric projection
onto a hyperbolic hyperplane in the hyperboloid model. Recall that the
$n$-dimensional hyperbolic space is given by
\begin{equation*}
\mathbb{H}^n
:=
\left\{
p\in\mathbb{R}^{n+1}:
\langle\!\!\langle p,p\rangle\!\!\rangle=-1,\;
p_{n+1}>0
\right\},
\end{equation*}
where
$
\langle\!\!\langle x,y\rangle\!\!\rangle
:=
x^\top Jy$,
$
J:=\operatorname{diag}(1,\ldots,1,-1),
$
denotes the Lorentzian inner product. The tangent space at
\(p\in\mathbb{H}^n\) is
$
T_p\mathbb{H}^n
=
\left\{
w\in\mathbb{R}^{n+1}:
\langle\!\!\langle p,w\rangle\!\!\rangle=0
\right\}.
$
The Riemannian metric on \(\mathbb{H}^n\) is the restriction of the
Lorentzian inner product to the tangent spaces. Thus, 
$
\langle u,w\rangle_p
=
\langle\!\!\langle u,w\rangle\!\!\rangle
$ holds for all \(u,w\in T_p\mathbb{H}^n\).
For simplicity, we omit the subscript whenever the base point is clear
from the context.
The Riemannian distance between \(p,q\in\mathbb{H}^n\) is given by
$d(p,q)
=
\operatorname{arcosh}
\bigl(
-\langle\!\!\langle p,q\rangle\!\!\rangle
\bigr).$
Moreover, the exponential and logarithmic maps at
\(q\in\mathbb{H}^n\) are given by
\[
\exp_qv
=
\begin{cases}
q,
& v=0,\\[3pt]
\cosh(\|v\|)q
+
\sinh(\|v\|)
\dfrac{v}{\|v\|},
& v\neq0,
\end{cases}
\]
and
\[
\log_qp
=
\begin{cases}
0,
& p=q,\\[3pt]
\dfrac{d(q,p)}
{\sqrt{\langle\!\!\langle p,q\rangle\!\!\rangle^2-1}}
\bigl(
p+
\langle\!\!\langle p,q\rangle\!\!\rangle q
\bigr),
& p\neq q.
\end{cases}
\]
Let \(q\in\mathbb{H}^n\) and \(v\in T_q\mathbb{H}^n\setminus\{0\}\).
We consider the set
$\mathcal{H}_{q,v}
:=
\left\{
p\in\mathbb{H}^n:
\langle v,\log_q p\rangle=0
\right\}.$
For \(p\neq q\), using the expression of the logarithmic map and the
fact that
\(\langle\!\!\langle v,q\rangle\!\!\rangle=0\), we obtain
$$
\langle v,\log_q p\rangle
=
\frac{d(q,p)}
{\sqrt{\langle\!\!\langle p,q\rangle\!\!\rangle^2-1}}
\,
\left\langle\!\!\left\langle
v,
p+\langle\!\!\langle p,q\rangle\!\!\rangle q
\right\rangle\!\!\right\rangle\\
=
\frac{d(q,p)}
{\sqrt{\langle\!\!\langle p,q\rangle\!\!\rangle^2-1}}
\,
\langle\!\!\langle v,p\rangle\!\!\rangle.
$$
Since \(q\in\mathcal{H}_{q,v}\), the above identity also yields
$
\mathcal{H}_{q,v}
=
\left\{
p\in\mathbb{H}^n:
\langle\!\!\langle v,p\rangle\!\!\rangle=0
\right\}.
$
Hence, \(\mathcal{H}_{q,v}\) is the intersection of \(\mathbb{H}^n\)
with the Lorentzian hyperplane \(v^\perp_h := \left\{ x \in \mathbb{R}^{n+1} : \langle\!\!\langle v, x \rangle\!\!\rangle = 0 \right\}\). In particular,
\(\mathcal{H}_{q,v}\) is a totally geodesic hyperplane of
\(\mathbb{H}^n\).

We next derive an explicit expression for the metric projection onto
\(\mathcal{H}_{q,v}\).

\begin{lemma}
\label{lem:projection-hyperplane}
Let \(q\in\mathbb{H}^n\) and
\(v\in T_q\mathbb{H}^n\setminus\{0\}\). Then, for every
\(p\in\mathbb{H}^n\), the metric projection of \(p\) onto
\(\mathcal{H}_{q,v}\) is given by
\[
\mathcal{P}_{\mathcal{H}_{q,v}}(p)
=
\frac{\phi(p)}
{\sqrt{-\langle\!\!\langle\phi(p),\phi(p)\rangle\!\!\rangle}},
\;\; 
\phi(p)
:=
p
-
\frac{\langle\!\!\langle v,p\rangle\!\!\rangle}
{\langle\!\!\langle v,v\rangle\!\!\rangle}
v.
\]
\end{lemma}

\begin{proof}
Since \(v\in T_q\mathbb{H}^n\setminus\{0\}\), we have
$
\langle\!\!\langle v,v\rangle\!\!\rangle>0.
$
By definition of $\phi$, we have
$
\langle\!\!\langle v,\phi(p)\rangle\!\!\rangle
=
0.
$
Thus, \(\phi(p) \in v_h^\perp\). Moreover, \(\langle\!\!\langle \phi(p), \phi(p) \rangle\!\!\rangle = \langle\!\!\langle p,p \rangle\!\!\rangle - \langle\!\!\langle v,p \rangle\!\!\rangle^2 / \langle\!\!\langle v,v \rangle\!\!\rangle = -1 - \langle\!\!\langle v,p \rangle\!\!\rangle^2 / \langle\!\!\langle v,v \rangle\!\!\rangle < 0\).
Consequently,
\[
\widehat p
:=
\frac{\phi(p)}
{\sqrt{-\langle\!\!\langle\phi(p),\phi(p)\rangle\!\!\rangle}}
\]
satisfies
$
\langle\!\!\langle\widehat p,\widehat p\rangle\!\!\rangle=-1
$
and
$
\langle\!\!\langle v,\widehat p\rangle\!\!\rangle=0.
$
Choosing the normalization with positive last coordinate, we obtain
$
\widehat p\in\mathcal{H}_{q,v}.
$
It remains to show that \(\widehat p\) minimizes the distance from \(p\)
over \(\mathcal{H}_{q,v}\). Let \(z\in\mathcal{H}_{q,v}\). Since
\(z\in v^\perp_h\), we have
$
\langle\!\!\langle p,z\rangle\!\!\rangle
=
\langle\!\!\langle\phi(p),z\rangle\!\!\rangle.
$
Since
$
\phi(p)
=
\sqrt{-\langle\!\!\langle\phi(p),\phi(p)\rangle\!\!\rangle}
\,\widehat p,
$
it follows that
\[
-\langle\!\!\langle p,z\rangle\!\!\rangle
=
\sqrt{-\langle\!\!\langle\phi(p),\phi(p)\rangle\!\!\rangle}
\,
\bigl(
-\langle\!\!\langle\widehat p,z\rangle\!\!\rangle
\bigr).
\]
Since \(\widehat p,z\in\mathbb{H}^n\),
$
-\langle\!\!\langle\widehat p,z\rangle\!\!\rangle
\geq1,
$
with equality if and only if \(z=\widehat p\). Therefore,
$
-\langle\!\!\langle p,z\rangle\!\!\rangle
\geq
-\langle\!\!\langle p,\widehat p\rangle\!\!\rangle.
$
Since \(\operatorname{arcosh}\) is strictly increasing,
$
d(p,z)
\geq
d(p,\widehat p)
$
holds for all $ z\in\mathcal{H}_{q,v}.
$
Thus,
$
\mathcal{P}_{\mathcal{H}_{q,v}}(p)
=
\widehat p,
$
which proves the result.
\end{proof}

\section{Sufficient Conditions for Assumption~\ref{ass:conv.su.}}\label{App:ass:conv.su.}
The purpose of this section is to provide conditions under which Assumption \ref{ass:conv.su.} is satisfied. We first characterize the convexity of the sets $\mathcal S_{q,v}$ in terms of their boundary  given by: 
\begin{equation}\label{eq:geodesic-hyperplane-boundary}
\mathcal H_{q,v}:=\mathcal S_{q,v}\cap\mathcal S_{q,-v}=\left\{p\in\mathbb M:\left\langle v,\log_qp\right\rangle=0\right\}=\exp_q(v^\perp),
\end{equation}
where $v^\perp:=\left\{\xi\in T_q\mathbb M:\langle v,\xi\rangle=0\right\}$.
The next proposition gives the precise relation between the convexity of
the two half-spaces and that of their common boundary.

\begin{proposition}
\label{prop:boundary-characterization}
Let 
$\dim\mathbb M\geq 2$, let $q\in\mathbb M$, and let
$v\in T_q\mathbb M\setminus\{0\}$. Then the following statements are
equivalent:
\begin{enumerate}
\item[(i)] $\mathcal H_{q,v}$ is  convex; 
\item[(ii)] both $\mathcal S_{q,v}$ and $\mathcal S_{q,-v}$ are convex.
\end{enumerate}
\end{proposition}
\begin{proof}
By the Cartan--Hadamard theorem, $\exp_q:T_q\mathbb M\to\mathbb M$ is a diffeomorphism. Hence $\mathcal H_{q,v}=\exp_q(v^\perp)$ is a closed embedded hypersurface. It separates $\mathbb M$ into the two open connected sets
$$
\mathcal{H}_{q,v}^-:=\left\{p\in\mathbb M:\left\langle v,\log_qp\right\rangle<0\right\} \;\; {\rm and }\;\; \mathcal{H}_{q,v}^+:=\left\{p\in\mathbb M:\left\langle v, \log_qp\right\rangle>0\right\}.
$$
In particular, $\mathcal S_{q,v}=\mathcal{H}_{q,v}^-\cup\mathcal H_{q,v}$ and $\mathcal S_{q,-v}=\mathcal{H}_{q,v}^+\cup\mathcal H_{q,v}$.
Assume first that $\mathcal H_{q,v}$ is convex. Let $p,r\in\mathcal S_{q,v}$, and let $\gamma:[0,1]\to\mathbb M$ be the unique geodesic segment joining $p$ to $r$. Suppose, by contradiction, that $\gamma([0,1])$ is not contained in $\mathcal S_{q,v}$. Then $I:=\left\{ t\in[0,1]:\gamma(t)\in\mathcal{H}_{q,v}^+ \right\}$ is nonempty. Let $(a,b)$ be a connected component of $I$. Since $\gamma(0),\gamma(1)\in\mathcal S_{q,v}$, continuity gives $\gamma(a),\gamma(b)\in\mathcal H_{q,v}$ and $\gamma(t)\in\mathcal H_{q,v}^+$, for every $a<t<b$. The convexity of $\mathcal H_{q,v}$ implies that the geodesic segment joining $\gamma(a)$ to $\gamma(b)$ lies in $\mathcal H_{q,v}$. Since geodesic segments joining pair of points are unique on a Hadamard manifold, this segment is the restriction of $\gamma$ to $[a,b]$. Therefore, $\gamma([a,b])\subset\mathcal H_{q,v}$, which contradicts $\gamma((a,b))\subset\mathcal{H}_{q,v}^+$. Thus, $\mathcal S_{q,v}$ is convex. By the same argument we have that $\mathcal S_{q,-v}$ is also convex. 
Conversely, if both sets $\mathcal S_{q,v}$ and $\mathcal S_{q,-v}$ are convex, then their intersection $\mathcal H_{q,v}= \mathcal S_{q,v}\cap\mathcal S_{q,-v}$ is  convex. 
\end{proof}

The above boundary characterization immediately shows that no curvature constancy assumption is needed in dimension two.

\begin{proposition} \label{prop:hadamard-surface}
Let $\mathbb M$ be a Hadamard surface. Then, for every $q\in\mathbb M$ and every $v\in T_q\mathbb M\setminus \{0\}$, the set $\mathcal S_{q,v}$ is convex. 
\end{proposition}

\begin{proof}
Since $\dim\mathbb M=2$, the space $v^\perp$ is one-dimensional. Choose a unit vector $w\in T_q\mathbb M$ such that $v^\perp=\operatorname{span}\{w\}$. Then $\mathcal H_{q,v}=\exp_q(v^\perp)=\left\{\exp_q(tw):t\in\mathbb R\right\}$ is the image of a complete geodesic. 
If two points of $\mathcal H_{q,v}$ are given, the subarc of this geodesic joining them is a geodesic segment in $\mathbb M$. Since an ambient geodesic segment between two points of a Hadamard manifold is unique, the previous subarc is the unique geodesic segment joining those points. Hence, $\mathcal H_{q,v}$ is  convex and the result follows from Proposition~\ref{prop:boundary-characterization}.
\end{proof}

Next, using the characterization established above, we show that Assumption~\ref{ass:conv.su.} is also satisfied on Hadamard manifolds of constant sectional curvature and dimension at least three.

\begin{proposition}
\label{lem:hyperp.conv}
Let $\mathbb M$ be an $n$-dimensional Hadamard manifold with constant sectional curvature, where $n\geq3$. Then, for every $q\in\mathbb M$ and every $v\in T_q\mathbb M\setminus \{0\}$, the set $\mathcal S_{q,v}$ is  convex.
\end{proposition}

\begin{proof}
Set $V:=v^\perp\subset T_q\mathbb M$. Since $\dim\mathbb M=n\geq3$, we have $2\leq\dim V<\dim\mathbb M$. Since $\mathbb M$ is complete and has constant sectional curvature, by
\cite[Theorem~2.2.4]{BerndtSuh2022}, there exists a complete totally geodesic submanifold $N$ of $\mathbb M$ such that $q\in N$ and $T_qN=V$.  

We claim that $N=\exp_q(V)=\mathcal H_{q,v}$. Indeed, since $N$ is complete,   for $u\in V=T_qN$ the intrinsic geodesic $\eta_u$ of $N$ satisfying $\eta_u(0)=q$ and  $\eta'_u(0)=u$ is defined for all time. Since $N$ is totally geodesic, $\eta_u$ is also a geodesic of $\mathbb M$. By uniqueness of the geodesic with prescribed initial data, $\eta_u(t)=\exp_q(tu)$. In particular, $\exp_q(u)\in N$, and hence $\exp_q(V)\subset N$.  Conversely, let $p\in N$. Since $N$ is complete, the Hopf--Rinow theorem yields an intrinsic minimizing geodesic $\alpha:[0,1]\to N$ joining $q$ to $p$. Setting $w:=\alpha'(0)\in T_qN=V$, and using once again the fact that $N$ is totally geodesic, we obtain $p=\alpha(1)=\exp_q(w)$. Thus, $N\subset\exp_q(V)$. Consequently, $N=\exp_q(V)=\mathcal H_{q,v}.$

It remains to prove that $N=\mathcal H_{q,v}$ is geodesically convex. Let $p,r\in N$. Since $N$ is complete, the Hopf--Rinow theorem provides an intrinsic minimizing geodesic in $N$ joining $p$ to $r$. Again, by the fact that $N$ is totally geodesic, this curve is also an ambient geodesic of $\mathbb M$. Since $\mathbb M$ is a Hadamard manifold, the ambient geodesic segment joining $p$ to $r$ is unique. It must therefore coincide with the geodesic contained in $N$. Hence, the unique ambient geodesic segment joining any two points of $N$ remains in $N$, thus $N$ is  convex.

Since $\mathcal H_{q,v}$ is geodesically convex, the Proposition~\ref{prop:boundary-characterization} implies that $\mathcal S_{q,v}$ is  convex. 
\end{proof}

\begin{remark}
\label{rem:cartan}
Suppose that $n:=\dim\mathbb{M}\geq3$ and that $\mathcal{S}_{q,v}$ is convex for \emph{every} $q\in\mathbb{M}$ and \emph{every} $v\in T_q\mathbb{M}\setminus\{0\}$. Applying the hypothesis to $v$ and $-v$, Proposition~\ref{prop:boundary-characterization} shows that every hypersurface $\mathcal{H}_{q,v}=\exp_q(v^\perp)$ is  convex. A convex closed embedded hypersurface is totally geodesic. Indeed the unique ambient geodesic segment joining two of its points lies in it, so its intrinsic geodesics are, locally, ambient geodesics. Since every hyperplane $W\subset T_q\mathbb{M}$ is of the form $v^\perp$, we conclude that through every point of $\mathbb{M}$ and tangent to every hyperplane there passes a totally geodesic hypersurface. This is Cartan's \emph{axiom of $(n-1)$-planes}, and a classical theorem of Cartan asserts that a Riemannian manifold of dimension $n\geq3$ satisfying the axiom of $r$-planes for some $2\leq r\leq n-1$ has constant sectional curvature, see \cite[Theorem~1.8]{Dajczer1990} . Hence, in dimension $n\geq3$, the uniform version of Assumption~\ref{ass:conv.su.} holds only on the space forms of Proposition~\ref{lem:hyperp.conv}. This explains why Assumption~\ref{ass:conv.su.} is formulated only along the pairs $(q_k,v_k)$ generated by the algorithm: on a manifold of nonconstant curvature it can hold for the particular directions produced by a given objective function although it fails for generic directions. In dimension two no such restriction exists (Proposition~\ref{prop:hadamard-surface}).
\end{remark}

\section*{Declarations}

\noindent
Conflict of interest: The authors declare that they have no conflict of interest.

\medskip
\noindent
Ethics approval and consent to participate: Not applicable.

\medskip
\noindent
Consent for publication: Not applicable.

\medskip
\noindent
Materials availability: Not applicable.

\medskip
\noindent
Code availability: Not applicable.

\bibliographystyle{plain}

\bibliography{references}

\end{document}